\documentstyle[12pt]{article} 
\begin{document} 
\newtheorem{Def}{Definition}[section] 
\newtheorem{thm}{Theorem}[section] 
\newtheorem{lem}{Lemma}[section] 
\newtheorem{rem}{Remark}[section] 
\newtheorem{prop}{Proposition}[section] 
\newtheorem{cor}{Corollary}[section] 
\newtheorem{conj}{Conjecture}[section]
\newtheorem{question}{Question}[section]
\title 
{Conformally invariant fully nonlinear
elliptic equations and isolated singularities }
\author{YanYan Li\thanks{Partially
 supported by
       NSF grant DMS-0401118.}\\
       Department of Mathematics\\
       Beijing Normal University\\
       Beijing 100875\\
       China\\
       and\\
       \\
       Department of Mathematics\\
       Rutgers University\\
       110 Frelinghuysen Road\\
       Piscataway, NJ 08854\\
       USA
       }
\date{} 
\maketitle 
\input { amssym.def} 

\setcounter{section}{0}

\section{ Introduction} 

There has been much work on conformally invariant
fully nonlinear elliptic equations
and  applications to geometry and topology.
See for instance \cite{V1},
 \cite{CGY}, \cite{C}, \cite{LL1},  \cite{LL3},
 \cite{GV}, and the references therein.
 In this and a companion paper \cite{Li1}
 we address some analytical issues concerning 
 these equations.

 For $n\ge 3$, consider
 \begin{equation}
 -\Delta u=n(n-2) u^{ \frac{n+2}{n-2}},
 \qquad \mbox{on} \quad \Bbb {R}^n.
 \label{eq1new}
 \end{equation}
 The method of moving planes was used by  Gidas, Ni and Nirenberg (\cite{GNN})
 in proving that
  any  positive $C^2$ solution  of
  (\ref{eq1new}) satisfying
  $\int_{\Bbb R^n}u^{ \frac {2n}{n-2} }<\infty$
  must be of the form
  $$
  u(x)=
  \left(\frac a {1+a^2|x-\bar{x}|^2}\right)^{\frac{n-2}{2}},
  $$
  where $a>0$ and $\bar x\in \Bbb R^n$.
  The hypothesis $\int_{\Bbb R^n}u^{ \frac {2n}{n-2} }<\infty$
  was removed
  by Caffarelli, Gidas and Spruck
  (\cite{CGS});
   this is important
  for applications.
This latter result
was extended to general
 conformally invariant
fully nonlinear 
second order 
elliptic equations
 in joint work with Aobing Li
\cite{LL3}, see also \cite{LL2} and \cite{LL4}.
For earlier results on the Liouville type theorems, see
 \cite{LL4} for a description.
Behavior near the origin of positive solutions of
$\displaystyle{-\Delta u=u^{ \frac{n+2}{n-2}}}$ 
in a punctured ball is also analyzed in \cite{CGS}.
Among other things,
we extend in this paper  a number of results
in \cite{CGS} to  general
 conformally invariant
 second order 
 fully nonlinear elliptic equations.
New techniques are 
 developed in the present paper.  Some of these,
in particular Theorem \ref{thm10}, 
have been used  in the companion paper \cite{Li1} to 
study  general degenerate 
 conformally invariant
 fully nonlinear elliptic equations.

Let ${\cal S}^{n\times n}$ denote the set of
$n\times n$ real symmetric matrices,
 ${\cal S}^{n\times n}_+$ denote the subset of ${\cal S}^{n\times n}$
 consisting of positive definite matrices,
  $O(n)$ denote
  the set of $n\times n$ real orthogonal matrices,
   $U\subset {\cal S}^{n\times n}$
   be  an open set,
   and  $F\in C^1(U)\cap C^0(\overline U)$.
   
      We list below a number of properties of $(F, U)$.
Subsets of these properties
         are used in various lemmas, propositions  and theorems:
\begin{equation}
O^{-1}UO=U,\qquad
\forall\ O\in O(n),
\label{21-1}
\end{equation}
\begin{equation}
U\cap \{M+tN\ |\
0<t<\infty\}\
\mbox{is convex}, \
\forall\ M\in {\cal S}^{n\times n},
N\in {\cal S}^{n\times n}_+,
\label{21-2}
\end{equation}
\begin{equation}
M\in U\ \mbox{and}\ N\in {\cal S}^{n\times n}_+\
\mbox{implies}\ M+N\in U,
\label{M1-0}
\end{equation}
\begin{equation}
M\in U\ \mbox{and}\ a>0\
\mbox{implies}\ aM\in U,
\label{M1-0new}
\end{equation}
\begin{equation}
F(O^{-1}MO)=F(M),\quad
\forall\ M\in U, \forall\ O\in O(n),
\label{21-3}
\end{equation}
\begin{equation}
(F_{ij}(M))>0,\qquad
\forall\ M\in U,
\label{21-4}
\end{equation}
where $F_{ij}(M):=\frac{\partial F}{ \partial M_{ij} }(M)$, and,
for some $\delta>0$,
\begin{equation}
F(M)\ne 1\ \ \forall\
M\in U\cap
\{M\in {\cal S}^{n\times n}\
|\
\|M\|:= (\sum_{i,j}M_{ij}^2)^{\frac 12}<\delta\}.
\label{21-5}
\end{equation}

Examples of such $(F,U)$ include those given by the
elementary symmetric functions.
 For  $1\le k\le n$, let
  $$
   \sigma_k(\lambda)=\sum_{1\le i_1<\cdots <i_k\le n}\lambda_{i_1}
    \cdots \lambda_{i_k}
     $$
      be the $k-$th elementary symmetric function and let
       $
        \Gamma_k$ be the connected component
	 of $\{\lambda\in \Bbb R^n\ |\
	  \sigma_k(\lambda)>0\}$ containing the positive cone $\Gamma_n
	   :=\{\lambda=(\lambda_1, \cdots, \lambda_n)\ |\
	    \lambda_i>0\}$.
	    Let
	    $$
	    U_k:=\{ M\in {\cal S}^{n\times n}\
	    |\ \lambda(M)\in \Gamma_k\},
	    $$
	    and
	    $$
	    F_k(M):=\sigma_k(\lambda(M))^{ \frac 1k},
	    $$
	    where $\lambda(M)$ denotes the eigenvalues of $M$.
	    Then $(F,U)=(F_k, U_k)$ satisfy all the above listed
	    properties, see for instance \cite{CNS}.
Taking $k=1$, equation
$$
F_1(A^u)=1
$$
amounts to, modulo a harmless positive constant,
$$
-\Delta u= u^{ \frac {n+2}{n-2} }.
$$
 Here and throughout the paper we use notation
  $$
   A^u= -\frac{2}{n-2}u^{  -\frac {n+2}{n-2} }
    \nabla^2u+ \frac{2n}{(n-2)^2}u^ { -\frac {2n}{n-2} }
     \nabla u\otimes\nabla u-\frac{2}{(n-2)^2} u^ { -\frac {2n}{n-2} }
      |\nabla u|^2I,
       $$
       where $\nabla u$ denotes the gradient of $u$ and
       $\nabla^2 u$ denotes the Hessian of $u$.

Other, much more general, examples are as follows.
Let
$$
\Gamma\subset \Bbb R^n \ \mbox{be
an open convex symmetric cone  with vertex at the origin }
$$
satisfying
$$
\Gamma_n
\subset \Gamma\subset \Gamma_1:=\{\lambda\in\Bbb{R}^n|\sum\limits_i
\lambda_i>0\}.
$$
Naturally, $\Gamma$ being symmetric
means $(\lambda_1, \lambda_2,\cdots, \lambda_n)\in
\Gamma$ implies $(\lambda_{i_1}, \lambda_{i_2}, \cdots , \lambda_{i_n})
\in \Gamma$ for any permutation
 $(i_1, i_2, \cdots , i_n)$ of
  $(1,2,\cdots, n)$.

  Let
  $$
  f\in C^1(\Gamma)\cap
   C^0(\overline \Gamma)
    $$
     satisfy
      $$
       f|_{\partial\Gamma}=0,\qquad
        \nabla f\in\Gamma_n\ \mbox{on}\ \Gamma,
	 $$
	  and
	  $$
	   f(s\lambda)=sf(\lambda),\qquad \forall\ s>0\
	    \mbox{and}\
	      \lambda\in \Gamma.
	       $$

	       With such $(f, \Gamma)$,  let
	       $$
	       U:=\{M\in {\cal S}^{n\times n}\ |\
	       \lambda(M)\in \Gamma\},
	       $$
	       and
	       $$
	       F(M):=f(\lambda(M)).
	       $$
	       Then $(F, U)$ satisfies all the above listed properties.
	       In fact, for all these $(F,U)$, $A^u\in U$ implies
	       $\Delta u\le 0$.
	       So for these $(F,U)$, the assumption
	       $\Delta u\le 0$ in various theorems in this paper is automatically
	       satisfied.
	       We note that in all these examples, $F$ is actually concave in $U$,
	       but this property is not needed in results in this paper.

Throughout the paper we use $B_a(x)\subset \Bbb R^n$ to denote 
the ball of radius $a$ and centered at $x$,
and $B_a=B_a(0)$.
Also, unless otherwise stated, the dimension $n$ is bigger than $2$.

\begin{thm}
Let $U\subset {\cal S}^{n\times n}$ be an open set satisfying (\ref{21-1})
and  (\ref{21-2}), let $F\in C^1(U)$ satisfy (\ref{21-3}), (\ref{21-4}) and
(\ref{21-5}).  Assume that $u\in C^2(\Bbb R^n\setminus B_{\frac 12})$
satisfy
\begin{equation}
u>0, \ \ \Delta u\le 0,
\qquad\mbox{in}\ \Bbb R^n\setminus B_{\frac 12},
\label{22-1}
\end{equation}
and
\begin{equation}
F(A^u)=1, \ \ A^u\in U, 
\qquad\mbox{in}\ \Bbb R^n\setminus B_{\frac 12}.
\label{22-3}
\end{equation}
Then
\begin{equation} 
\limsup_{ |x|\to \infty}
|x|^{ \frac {n-2}2 }u(x)<\infty.
\label{22-4}
\end{equation}
\label{thm2}
\end{thm}

\begin{rem} 
For $(F,U)=(F_1, U_1)$, (\ref{22-4})
was proved in \cite{CGS}.
\label{rem4}
\end{rem}

\begin{rem} Gonzalez  in \cite{d2} and Han in \cite{H}
 studied for certain  $(F_k, U_k)$
solutions with isolated singularities
which have finite volume, and Gonzalez  in \cite{d1} studied
subcritical  $(F_k, U_k)$                                            
solutions with isolated singularities. Chang, Han and Yang studied
in \cite{CHY} radial solutions on  annular domains including
punctured balls and $\Bbb R^n$.
See these papers for precise statements and details.
\end{rem}

\begin{rem} If $diag(-\frac 12, \frac 12, \cdots, \frac 12)\in
U$, then the upper bound  (\ref{22-4})  is sharp
in the sense that the exponent $\frac {n-2}2$ can not be larger.
This is because
$$
\lambda(A^u)\equiv \left\{
-\frac 12, \frac 12, \cdots, \frac 12\right\},
\qquad \mbox{on}\ \Bbb R^n\setminus\{0\}
$$
for $u(x)=|x|^{ \frac {2-n}2 }$. In particular, (\ref{22-4}) is sharp
for $(F,U)=(F_k, U_k)$ for $1\le k<\frac n2$.
See Section 8 for details.
\label{rem1.3}
\end{rem}

\begin{rem} Condition (\ref{21-5}) can not be dropped since $u\equiv
constant$ could be a solution.
\end{rem}
\begin{rem}
Instead of (\ref{22-4}),  what we have actually  proved is
$$
\sup_{|x|\ge 1} |x|^{ \frac{n-2}2}
u(x)\le C,
$$
for some $C$ explicitly given in terms of $
\displaystyle{
\min_{\partial B_1}u}$ and $
n$.  This can be seen from the proof 
of Theorem \ref{thm2}.
\label{rem32-1}
\end{rem}

Replacing $u$ by $|x|^{2-n}u(\frac x{|x|^2})$ and using the 
conformal invariance property of $F(A^u)$
--- see for  example line 9 on page
1431 of \cite{LL1}, it is easy to see that
Theorem \ref{thm2} is equivalent to

\noindent{\bf Theorem \ref{thm2}$'$}.\
{\it
Let $U\subset {\cal S}^{n\times n}$ be an open set satisfying (\ref{21-1})
and  (\ref{21-2}), let $F\in C^1(U)$ satisfy (\ref{21-3}), (\ref{21-4}) and
(\ref{21-5}).  Assume that $u\in C^2(B_2\setminus\{0\})$
satisfy
}
\begin{equation}
u>0, \ \Delta u\le 0, 
\qquad\mbox{in}\ B_2\setminus\{0\},
\label{22-1prime}
\end{equation}
{\it and }
\begin{equation}
F(A^u)=1, \ \ A^u\in U,
\qquad\mbox{in}\  B_2\setminus\{0\}.
\label{22-4prime}
\end{equation}
{\it Then }
$$
\limsup_{ |y|\to 0}
|y|^{ \frac {n-2}2 }u(y)<\infty.
$$

\bigskip

\begin{thm} Let $U\subset {\cal S}^{ n\times n}$ be an open set satisfying
(\ref{21-1}) and (\ref{21-2}), and let
$F\in C^1(U)$ satisfy (\ref{21-3}) and (\ref{21-4}).
Assume that $u\in C^2(\Bbb R^n\setminus\{0\})$ satisfy
$$
u>0, \ \ \Delta u\le 0, \qquad\mbox{in}\
\Bbb R^n\setminus\{0\},
$$
$$
F(A^u)=1, \quad A^u\in U, \quad \mbox{in}\ \Bbb R^n\setminus\{0\},
$$
and
\begin{equation}
u\ \mbox{can not be extended as a }\ C^2 \ \mbox{positive
function satisfying
$A^u\in U$ near the origin}.
\label{34-1}\end{equation}
Then $u$ is radially symmetric about the origin.
\label{thm3}
\end{thm}

\begin{rem} For $(F,U)=(F_1, U_1)$, the result was proved in \cite{CGS}.
\end{rem}

\begin{thm}
Let $U\subset {\cal S}^{n\times n}$ be an open set satisfying (\ref{21-1})
and  (\ref{21-2}), and let $F\in C^1(U)$ satisfy (\ref{21-3}), (\ref{21-4}) and
(\ref{21-5}).  Assume that $u\in C^2(B_2\setminus\{0\})$
satisfies (\ref{22-1prime}) and (\ref{22-4prime}). Then, for some
constant $\epsilon>0$,
\begin{equation}
u_{x,\lambda}(y)\le u(y),\qquad
\forall\ 0<\lambda< |x|\le \epsilon,
\ |y-x|\ge \lambda, 0<|y|\le 1.
\label{radial}
\end{equation}
Consequently, 
 for some positive constant $C$,
\begin{equation}
\bigg|\frac {u(x)}{ u(y)}-1\bigg|
\le Cr,\qquad \forall\ 0<r=|x|=|y|< 1.
\label{vv9}
\end{equation}
\label{thm11}
\end{thm}

\begin{rem} For $(F,U)=(F_1, U_1)$, the result was proved in \cite{CGS}.
\end{rem}

\begin{rem} In view of Remark \ref{rem32-1}, we can  obtain 
explicit dependence of $\epsilon$ and $C$ in terms of $
\displaystyle{
\min_{\partial B_1}u
}$ and $
n$. 
    With such explicit dependence,
Theorem \ref{thm3}
 follows from Theorem 
\ref{thm11} by rescaling a large ball to $B_2$ and
then sending the radius of the large ball to infinity.  
In doing this, the minimum of $\partial B_1$
of the rescaled function
is under control due to the fact
$
\displaystyle{
\liminf_{ |y|\to \infty} |y|^{n-2}u(y)>0.
}
$
We leave the details to interested readers.
\end{rem}

\begin{thm} Let
\begin{equation}
U\subset  {\cal S}^{n\times n}_+
\label{41-1}
\end{equation}
be an open set satisfying (\ref{21-1}) and (\ref{21-2}),
let $F\in C^1(U)$ satisfy (\ref{21-3}) and (\ref{21-4}),
and let $u\in C^2(B_2\setminus \{0\})$ satisfy
$$
u>0,\qquad \mbox{in}\ B_2\setminus\{0\},
$$
and
$$
F(A^u)=1, A^u\in U, \qquad \mbox{in}\
B_2\setminus\{0\}.
$$
Then $u$ can be extended as a positive Lipschitz function in
$B_1$.
\label{thm4}
\end{thm}

\begin{cor} The conclusion of Theorem \ref{thm4} holds
for $(F, U)=(F_n, U_n)$.
\label{cor1-1}
\end{cor}

\begin{thm} Let $U\subset {\cal S}^{n\times n}$ be an open 
set satisfying (\ref{M1-0}) and (\ref{M1-0new}).
We assume that there exists some
$\eta\in C^2(B_2\setminus\{0\})\cap C^0(B_2)$  satisfying
\begin{equation}
\eta(0)=0, \ \eta(x)>0\ \ \forall\ x\in B_2\setminus\{0\},
\label{M1-1}
\end{equation}
\begin{equation}
D^2\eta(x)\ \mbox{does not belong to}\
U,\qquad \forall\ x\in B_2\setminus\{0\}.
\label{M1-2}
\end{equation}
Suppose that $\xi\in C^0_{loc}(B_2\setminus\{0\})\cap
 L^\infty(B_2\setminus\{0\})$
satisfies 
\begin{equation}
\xi>0\ \mbox{in}\ B_2\setminus\{0\},
\ \ \Delta \xi\ge 0\
\mbox{in}\ B_2\setminus\{0\}
\ \mbox{in the distribution sense},
\label{M1-3}
\end{equation}
and there exist $\{\xi_i\}$ in $C^2(B_2\setminus\{0\})$ such that
\begin{equation}
\Delta \xi_i\ge 0\qquad\mbox{in}\
B_2\setminus\{0\},
\label{M2-0}
\end{equation}
\begin{equation}
D^2\xi_i\in \overline U,
\qquad \mbox{in}\ B_2\setminus\{0\},
\label{M2-1}
\end{equation}
\begin{equation}
\xi_i\to \xi,\qquad \mbox{in}\ C^0_{loc}(B_2\setminus\{0\}).
\label{M2-2}
\end{equation}
Then $\xi$ can be extended as a function in
$C^0(B_1)$ which satisfies 
\begin{equation}
\sup_{ B_1\setminus\{0\}}\xi\le \max_{\partial B_1}\xi,
\label{M2-3}
\end{equation}
\begin{equation}
|\xi(x)-\xi(y)|\le C(\eta)
\left[ \sup_{ B_1\setminus\{0\}}\xi
-\inf_{  B_1\setminus \{0\}  }\xi\right]
\left[ \eta(x-y)+\eta(y-x)\right],
\qquad
\forall\ x, y\in B_{\frac 14},
\label{M2-4}
\end{equation}
where $C(\eta)$ denotes some positive constant depending on
$\eta$.
\label{thmM1}
\end{thm}

\begin{cor}
For $B_2\subset \Bbb R^n$, $n\ge 1$,
let $k$ be an integer satisfying $\frac n2< k\le n$.  
We assume that $\xi\in C^2(B_2\setminus\{0\})\cap L^\infty(B_2\setminus\{0\})$
and
$$
\lambda(D^2\xi)\in \overline \Gamma_k\quad \mbox{in}\ B_2\setminus\{0\}.
$$
Then, for $\alpha =\frac{2k-n}k$, $\xi$ can be
extended as a function in $C^{0, \alpha}(B_1)$ and,
for any $0<a<2$, 
\begin{equation}
\|\xi\|_{  C^{0, \alpha}(B_a) }\le C(n,a)
\left(\sup_{ B_2\setminus\{0\} }  \xi-
\inf_{ B_2\setminus\{0\} }  \xi\right),
\label{nn88}
\end{equation}
where $C(n,a)$ is some positive constant depending only on $n$
and $a$.
\label{lemA-1}
\end{cor}

\begin{rem} Without the possible singularity
of $\xi$ at the origin,
(\ref{nn88}) was known, see theorem 2.7 in 
 \cite{TW} by Trudinger and Wang.
\end{rem}

\begin{cor} Let $U$ and $\eta$ be as in Theorem \ref{thmM1}.  Suppose
that $u\in C^0_{loc}(B_2\setminus\{0\})$ satisfies
$u>0$ in $B_2\setminus\{0\}$ and there exist $\{u_i\}$ in
$C^2(B_2\setminus\{0\})$,
$$
\Delta u_i\le 0,\ \ 
A^{u_i}\in \overline U,\qquad
\mbox{in}\ B_2\setminus \{0\},
$$
$$
u_i\to u\qquad \mbox{in}\
C^0_{loc}(B_2\setminus\{0\}).
$$
Then $\xi:=u^{ -\frac 2{n-2}}$ can
be extended as a function in $C^0(B_1)$ and
\begin{equation}
\sup_{ B_1\setminus\{0\}}\xi
\le \max_{ \partial B_1}\xi
=[\max_{ \partial B_1} u]^{  -\frac 2{n-2} },
\label{M12-1}
\end{equation}
\begin{equation}
|\xi(x)-\xi(y)|\le
C(\eta) [\min _{\partial B_1} u]^{ -\frac 2{n-2} }
[\eta(x-y)+\eta(y-x)],
\qquad\forall\ x,y\in B_{\frac 14}.
\label{M12-2}
\end{equation}
Consequently, either
\begin{equation}
0<\inf_{ B_1\setminus\{0\}}u
\le \sup_{ B_1\setminus\{0\}}u<\infty\ \ 
\mbox{and}\ u\in C^0(B_1),
\label{M13-1}
\end{equation}
or
\begin{equation}
\inf_{ x\in B_1\setminus\{0\}}[\eta(x)+\eta(-x)]^{ \frac {n-2}2 } u(x)>0.
\label{M13-2}
\end{equation}
\label{corM11}
\end{cor}

\begin{cor} Let $B_2\subset \Bbb R^n$ and let 
$k$ be   an integer satisfying $\frac n2 <k\le n$.  We assume
that $u\in C^2(B_2\setminus \{0\})$,
$
u>0$ and $\lambda(A^u)\in \overline \Gamma_k$
 on $B_2\setminus \{0\}$.  Then
  $\xi:= u^{ -\frac 2{n-2}}$ can be
   extended as a function in $C^{0, \alpha}(B_1)$,
    with $\alpha=
     \frac{2k-n} k\in (0, 1]$, and
    $$
       \|\xi\|_{  C^{0,
         \alpha}(B_{\frac 12})  }\le C(n)
	   [\min_{ \partial B_1} u]^{ -\frac 2{n-2} }.
$$
	        Consequently,
		   either
		     $$
0<\inf_{
B_1\setminus\{0\}} u
\le \sup_{       
B_1\setminus\{0\}} u
<\infty,
			   \ \ \mbox{and}\ u\in
			     C^{0,
			         \alpha}(B_1),
		$$
				     or
				       $$
		|x|^{ \frac {n-2}2 \alpha}u(x)
\ge \frac 1{C(n)}
[\min_{\partial B_1} u], 
\qquad \forall\ |x|<\frac 12.
				$$
								 \label{thm5}
								 \end{cor}

\begin{rem} 
The H\"older
regularity of $\xi$ was independently proved by
Gursky and Viaclovsky in \cite{GV}, which contains
some more general and  other
very nice results.
  Our proof
is 
 different.
\end{rem}

\begin{rem}
The H\"older exponent in Theorem \ref{thm5} is sharp, compare 
for instance results in \cite{CHY}.
\end{rem}

Our proofs of Theorem \ref{thm2},
Theorem \ref{thm3},  Theorem \ref{thm4} and
Theorem \ref{thmM1} make use of the following theorem and
 its
generalizations.

\begin{thm} 
Let $U\subset {\cal S}^{n\times n}$ be an open set satisfying (\ref{21-1})
and  (\ref{21-2}), and let $F\in C^1(U)$ satisfy (\ref{21-3})
and  (\ref{21-4}). 
 We assume that $u\in C^2(B_2\setminus\{0\})$ and
 $v\in C^2(B_2)$ satisfy
 $$
 u>v \qquad \mbox{in}\ B_2\setminus\{0\},
 $$
$$
 F(A^u)\ge 1, \quad A^u\in U,  \Delta u\le 0,
 \qquad  \mbox{in}\ B_2\setminus\{0\},
$$
$$
  F(A^v)\le  1, \quad A^v\in U, v>0, \qquad  \mbox{in}\ B_2.
$$
   Then
   \begin{equation}
   \liminf_{|x|\to 0}[u(x)-v(x)]>0.
   \label{conclusion}
   \end{equation}
\label{thm0}
\end{thm}

\begin{rem} As pointed out in \cite{LL3},
the  arguments
in \cite{LL1} together with Theorem \ref{thm0} yield
the Liouville type theorem in \cite{LL3}.
  The proof of the Liouville type theorem
in \cite{LL3} avoids such local result by using global information
of the entire solution $u$.  Our proof of  Theorem \ref{thm0} 
makes use of the crucial idea in the proof of the Liouville type
theorem in \cite{LL3} --- a delicate use of Lemma \ref{lemma00}.
\end{rem}

The conclusion of Theorem \ref{thm0} holds for elliptic operators 
with less invariance than the M\"obius group.
Let
 $T\in C^1(\Bbb R^+\times \Bbb R^n\times {\cal S}^{n\times n})$
satisfy
\begin{equation}
\left(-\frac{\partial T}{\partial u_{ij}}\right)>0
\qquad \mbox{on} \ \Bbb R_+\times \Bbb R^n\times {\cal S}^{n\times n},
\label{ba1}
\end{equation}
where $\Bbb R_+=(0, \infty)$.
With (\ref{ba1}), the operator $T(u, \nabla u, \nabla ^2 u)$ is elliptic.

For a positive function $v$, 
and for $x\in \Bbb R^n$ and $\lambda>0$, let
$$
v^{x,\lambda}(y):= \lambda^{ \frac {n-2}2 }
v(x+\lambda y), \qquad v^\lambda(y):=v^{0, \lambda}(y).
$$

We assume that the operator $T$ has the following invariance:
For any positive function $v\in C^2(\Bbb R^n)$ and for
any  $\lambda>0$,
\begin{equation}
T(v^{\lambda}, \nabla v^{\lambda}, \nabla^2 v^{\lambda})(\cdot)
\equiv T(v, \nabla v, \nabla^2 v)(\lambda \cdot)
\qquad \mbox{in}\ \Bbb R^n.
\label{1-1}
\end{equation}

\begin{rem} Let $\displaystyle{T(t,p,M):=S(t^{-\frac n{n-2}}p,
t^{ -\frac{n+2}{n-2} }M)}$ for some
$S\in C^0(\Bbb R^n\times {\cal S}^{n\times n})$.  Then
$T$ satisfies (\ref{1-1}).  See Lemma \ref{thmT}.
\end{rem}

\begin{thm} Let  $B_2\subset \Bbb R^n$
and let $T\in C^1(\Bbb R^+\times \Bbb R^n\times {\cal S}^{n\times n})$ 
satisfy
(\ref{1-1}).  We assume that $u\in C^2(B_2\setminus\{0\})$ and
$v\in C^2(B_2)$ satisfy
\begin{equation}
v>0\qquad \mbox{in}\ B_2,
\label{2-1}
\end{equation}
\begin{equation}
u>v \qquad \mbox{in}\ B_2\setminus\{0\},
\label{2-2}
\end{equation}
\begin{equation}
\Delta u\le 0  \qquad \mbox{in}\ B_2\setminus\{0\},
\label{2-3}
\end{equation}
\begin{equation}
T(u, \nabla u, \nabla^2 u)\ge 0\ge
T(v, \nabla v, \nabla^2 v)\qquad \mbox{in}\ B_2\setminus\{0\}.
\label{2-5}
\end{equation}
Then
\begin{equation}
\liminf_{|x|\to 0}[u(x)-v(x)]>0,
\label{3-2y}
\end{equation}
\label{thm1}
\end{thm}

\begin{rem} It is not difficult to see
from the proof of 
Theorem \ref{thm1} that we have only used
the following properties of $u, v$ and $T$:
$T\in C^1(\Bbb R^+\times \Bbb R^n\times {\cal S}^{n\times n})$, 
$u\in C^2(B_2\setminus\{0\})$ and
$v\in C^2(B_2)$
satisfy (\ref{2-1}), (\ref{2-2}),
(\ref{2-3}), and there exists
some $\epsilon_5>0$ such that
$$
T(u, \nabla u, \nabla^2u)\ge T(v^{x,\lambda},
\nabla v^{x,\lambda}, \nabla^2v^{x,\lambda})\quad
\mbox{on}\ B_{\epsilon_5}\setminus \{0\},\
\forall\ |x|<\epsilon_5, |\lambda-1|<\epsilon_5,
$$
and
for any
$|x|<\epsilon_5$, $|\lambda-1|<\epsilon_5$,
$|y|<\epsilon_5$ satisfying $u(y)=v^{x,\lambda}(y)$,
$\nabla u(y)=\nabla v^{x,\lambda}(y)$, $u\ge v^{x,\lambda}$
on $B_{ \epsilon_5}\setminus\{0\}$, we have
$$
\left(-\frac{\partial T}{\partial u_{ij}}
\bigg(u(y), \nabla u(y), \theta \nabla^2 u(y)+
(1-\theta)\nabla^2  v^{x,\lambda}(y)
\bigg)\right)>0\quad \forall\ 0\le \theta\le 1.
$$
\label{rem4-1}
\end{rem}

\begin{rem} Taking $F(A^u)-1$ as the operator $T$, 
the   properties in Remark \ref{rem4-1} are satisfied by the 
$u$ and $v$ in Theorem \ref{thm0} --- see arguments towards the end of the proof of lemma 2.1
in \cite{LL1}.  Therefore  Theorem \ref{thm0} is,
in view of Remark \ref{rem4-1},
a consequence of Theorem \ref{thm1}.
\label{rem4-2}
\end{rem}

The following follows from a classical  result in
\cite{E}:
 Let $E$ be a closed subset of $B_2$ of capacity
$0$ ---  the standard capacity with respect to the Dirichlet integral,
and 
let $u\in C^2(B_2\setminus E)$ and $v\in C^2(B_2)$
satisfy
$$
u>v\ \mbox{and}\ \Delta u\le 0\le \Delta v
\ \mbox{in}\ B_2\setminus E.
$$
Then
\begin{equation}
\liminf_{ dist(x, E)\to 0}[u(x)-v(x)]>0.
\label{hhh2}
\end{equation}

Theorem \ref{thm2} and Theorem \ref{thm1} can be
viewed as an extension of this
for $E=\{0\}$.
\begin{question} In Theorem \ref{thm1}, if we replace
$\{0\}$ by some $E$ with capacity $0$,  
does (\ref{hhh2}) still hold?  Maybe there is  a notion of
$T-$capacity for  (\ref{hhh2}) to hold for
zero $T-$capacity set $E$?
\label{question1}
\end{question}

A more concrete question is
\begin{question} Let $T$ be as in Theorem \ref{thm1}
or $F(A^u)$ be as in Theorem \ref{thm2}, and
let $E=E^k\subset B_2$ be an 
embedded closed smooth manifold of dimension $k$.
What is the $k^*(n,T)$ for which 
 (\ref{hhh2}) holds  for all $0\le k\le k^*$ ---
 with the hypotheses of Theorem \ref{thm1} or 
 Theorem \ref{thm2} for
 $\{0\}$ being changed in an obvious way 
 to that for $E^k$? What about
 for
 $$
 E^k=\{(x_1, \cdots, x_k, 0, \cdots, 0)\
 |\ \sum_{i=1}^n (x_i)^2=1\}?
 $$
\label{question2}
\end{question}

Another question is
\begin{question} For what classes of elliptic
operators $T(x, u, \nabla u, \nabla^2 u)$ the conclusion of 
Theorem \ref{thm1} holds?
\label{question3}
\end{question}

Concerning this question we will give in 
Corollary \ref{corA1} and Corollary \ref{cor9}
some operators with the property.

For a one variable function $\varphi$, we define,
instead of $v^{x,\lambda}$, 
$$
v^{x,\lambda}_\varphi(y)=
\varphi(\lambda)v(x+\lambda y), \qquad
v^\lambda_\varphi=v^{0,\lambda}_\varphi.
$$

\begin{thm}Let $\Omega\subset \Bbb R^n$ be a
bounded open set containing the origin $0$,
 $n\ge 2$, and 
let
$\varphi$ be a $C^1$ function
defined in a neighborhood of $1$
satisfying
$\varphi(1)=1$ and $\varphi'(1)>0$.
We assume that  $u\in C^0(\overline \Omega\setminus
\{0\})$, $v$ is $C^0$ in some open neighborhood of
$\overline \Omega$ and
$v$ is $C^1$ in a neighborhood of $0$,
\begin{equation}
v>0\qquad\mbox{in}\ \overline \Omega,
\label{2-1newnew}
\end{equation}
\begin{equation}
u>v \qquad \mbox{on}\ \overline \Omega\setminus\{0\},
\label{2-2newnew}
\end{equation}
\begin{equation}
\Delta u\le 0  \qquad \mbox{in}\ \Omega\setminus\{0\}.
\label{2-3newnew}
\end{equation}
Assume also that there exists some $\epsilon_3>0$ such that
 for any $|x|<\epsilon_3$ and $|\lambda-1|<\epsilon_3$,
\begin{equation}
\inf_{\Omega\setminus\{0\}}[
u-v^{x,\lambda}_\varphi]=0
\ \mbox{implies}\ \liminf_{ |y|\to 0}
[u-v^{x,\lambda}_\varphi](y)=0.
\label{gg5}
\end{equation}
Then (\ref{3-2y}) holds.
\label{thm7}
\end{thm}

\begin{thm} Under the hypotheses of Theorem \ref{thm7},
except changing $\varphi'(1)>0$ to 
$\varphi'(1)<0$.  Then (\ref{3-2y}) holds if
$\varphi'(1)<-1$.  If $-1\le \varphi'(1)<0$, 
either  (\ref{3-2y}) holds, or
 \begin{equation}
  \liminf_{|x|\to 0}[u(x)-v(x)]=0
   \label{3-2}
    \end{equation}
 and, 
for some $\epsilon>0$ and $V\in \Bbb R^n$,
\begin{equation}
\psi(v(x))
+V\cdot x\equiv 0, \qquad |x|<\epsilon,
\label{gg7}
\end{equation}
where
$$
\psi(s):= \int_{ v(0) }^s 
\frac {v(0)} t \varphi^{-1}\left(\frac {v(0)} t\right)dt.
$$
\label{thm8}
\end{thm}

We give a corollary which concerns
Question \ref{question3}.
Let 
$S\in C^1(\Bbb R^n\times {\cal S}^{n\times n})$ satisfy
\begin{equation}
\left( -\frac{\partial S}{ \partial M_{ij} }(p, M)\right)>0\qquad
\forall\ (p, M)\in \Bbb R^n\times {\cal S}^{ n\times n},
\label{AA1}
\end{equation}
and let, for $\beta\in \Bbb R\setminus\{0\}$,
\begin{equation}
T(t,p,M):=S\left( t^{  -\frac {1+\beta}\beta  }p,
t^{ -\frac{2+\beta}\beta }M\right),\qquad
(t,p,M)\in \Bbb R_+\times \Bbb R^n\times {\cal S}^{ n\times n}.
\label{AA2}
\end{equation}
\begin{cor} For $n\ge 2$, let
$S$, $\beta$ and $T$ be as above.  If $-1<\beta<0$, we further
require that
\begin{equation}
S(p, 0)\ge 0,\qquad \forall\ p\in \Bbb R^n.
\label{CCC2}
\end{equation}
Assume that $u\in C^2(B_2\setminus\{0\})$ and
$v\in C^2(B_2)$ satisfy (\ref{2-1}), (\ref{2-2}), (\ref{2-3}) and 
(\ref{2-5}).  Then (\ref{3-2y}) holds.
\label{corA1}
\end{cor}

Clearly, the arguments in the proofs of Theorem \ref{thm0},
 Theorem \ref{thm1}, Theorem \ref{thm7} and Theorem \ref{thm8}
 can be used to study some other problems.  For instance, let
$$
\Phi(v, x, \lambda; y):= \varphi(\lambda) v(x+
\xi(\lambda)y)+\psi(\lambda).
$$
We assume that $\varphi, \psi$ and $\xi$ are $C^1$ functions
near $1$ satisfying $\varphi(1)=\xi(1)=1$, $\psi(1)=0$, 
\begin{equation}
\varphi'\ge 0, \psi'\ge 0,
\varphi'+\psi'>0,\qquad \mbox{near}\ 1,
\label{Q1-1}
\end{equation}
and
\begin{equation}
\varphi'\xi+\varphi\xi'\ge 0,\qquad  \mbox{near}\ 1,
\label{Q1-2}
\end{equation}

Here is an extension of Theorem \ref{thm7}.
\begin{thm} Let $\varphi, \xi, \psi$ be as above, and
let $\Omega\subset \Bbb R^n$ be a bounded open
set containing the origin $0$, $n\ge 2$.  We assume that 
$u\in C^0(\overline \Omega\setminus\{0\})$,
$v$ is $C^0$ in some open neighborhood of
$\overline \Omega$ and $v$ is $C^1$ near the 
origin.  Assume also that (\ref{2-1newnew}), (\ref{2-2newnew})
and (\ref{2-3newnew}) hold, and
there exists some $\epsilon_4>0$ such that
for any $|x|<\epsilon_4$ and $|\lambda-1|<\epsilon_4$,
\begin{equation}
\inf_{\Omega\setminus\{0\}}[u-
\Phi(v, x, \lambda; \cdot)]=0
\ \mbox{implies}\
\liminf_{ |y|\to 0} [u(y)-\Phi(v, x, \lambda; y)]=0.
\label{gg5z}
\end{equation}
Then (\ref{3-2y}) holds.
\label{thm9}
\end{thm}

We now give some more 
operators $T$ for which the conclusion
of Theorem \ref{thm1} holds. For $S$ satisfying
(\ref{AA1}), we consider operators $T$ satisfying
one of the following.

\bigskip

(i)\ $T(t,p,M):=S(p,M)$.

\medskip

(ii)\ There exists $\epsilon>0$ such that
$$
sign\ T(\lambda t, \lambda p, \lambda M)=
sign\ T(t,p,M),\ \ 
\forall\ \ (\lambda, t, p, M)\in
(1-\epsilon, 1+\epsilon)\times \Bbb R_+\times \Bbb R^n\times {\cal S}^{n\times n}.
$$

\medskip

(iii)\ $\displaystyle{
T(t,p,M):=S\left( \frac 1{t+1}p, \frac 1{t+1}M\right),
\ (t,p,M)\in \Bbb R_+\times \Bbb R^n\times {\cal S}^{n\times n} }.$

\bigskip

\begin{cor} For $n\ge 2$, let $S\in C^1(\Bbb R^n\times {\cal S}^{n\times n})$
satisfy (\ref{AA1}) and let $T\in
C^1(\Bbb R_+\times \Bbb R^n\times {\cal S}^{n\times n})$
satisfy one of the above.  
Assume that $u\in C^2(B_2\setminus\{0\})$ and
$v\in C^2(B_2)$ satisfy (\ref{2-1}), (\ref{2-2}), (\ref{2-3}) and 
(\ref{2-5}).  Then (\ref{3-2y}) holds.
\label{cor9}
\end{cor}

Corollary \ref{cor9} follows from a more general

\begin{cor} For $n\ge 2$, let
 $T\in
C^1(\Bbb R_+\times \Bbb R^n\times {\cal S}^{n\times n})$
satisfy (\ref{ba1}), and let 
$u\in C^2(B_2\setminus\{0\})$ and
$v\in C^2(B_2)$ satisfy (\ref{2-1}), (\ref{2-2}), (\ref{2-3}) and 
(\ref{2-5}).  Assume that for some $\varphi$,
$\xi$, $\psi$ as in Theorem \ref{thm9} and
for some $\epsilon>0$,
\begin{equation}
T\left(\Phi(v,0,\lambda;\cdot),
\nabla \Phi(v,0,\lambda;\cdot), \nabla^2 \Phi(v,0,\lambda;\cdot)\right)
\le 0\quad\mbox{in}\ B_\epsilon\ \mbox{for all}\ |\lambda-1|<\epsilon.
\label{BB1}
\end{equation}
 Then (\ref{3-2y}) holds.
\label{cor1.6}
\end{cor}

The operators $T$ in Corollary \ref{cor9} satisfy 
the hypotheses of Corollary \ref{cor1.6}, see Section
9.

\medskip

In some applications, see \cite{Li1},
assumption (\ref{2-2newnew}) in Theorem \ref{thm9} needs to be weakened.
For this purpose, we give
\begin{thm} Let
 $\Omega\subset \Bbb R^n$ be a bounded open
set containing the origin $0$, $n\ge 2$.  We assume that
$u\in C^0(\Omega\setminus\{0\})$,
$v$ is $C^1$ in some open neighborhood of $\overline \Omega$,
 $v$ satisfies (\ref{2-1newnew}), $u$ satisfies
(\ref{2-3newnew}), and 
$$
u\ge v\ \mbox{in}\ \Omega\setminus\{0\}.
$$
  Assume also that  $\varphi, \xi, \psi$  are $C^1$ functions
  near $1$ satisfying $\varphi(1)=\xi(1)=1$, $\psi(1)=0$,
  $\varphi'(1)+\xi'(1)>0$, and
$$
 \varphi'(1)v(y)+\xi'(1)\nabla v(y)\cdot y
 +\psi'(1)>0,\qquad \forall\ y\in \overline \Omega,
$$
 and assume that 
  there exists some $\epsilon_4>0$ such that (\ref{gg5z}) holds 
  for any $|x|<\epsilon_4$ and $|\lambda-1|<\epsilon_4$.
  Then either (\ref{3-2y}) holds or
 $u= v= v(0)$ near the
  origin.
  \label{thm10}
  \end{thm}

As mentioned earlier, we make, as in \cite{LL3}, delicate use
of the following result.

\begin{lem} (\cite{LL3})\ 
For $n\ge 2$, $B_1\subset\Bbb {R}^n$, let $u\in L^1_{loc}(B_1\setminus
\{0\})$ be the solution of
\[
\Delta u\le 0\qquad\mbox{in}~B_1\setminus \{0\}
\]
in the distribution sense. 
Assume $\exists$ $a\in \Bbb R$ and $p\neq q\in\Bbb {R}^n$
such that
\[
u(x)\ge \max \{a+p\cdot x-\delta (x),a+q\cdot x-\delta
(x)\}\quad\forall x\in~B_1\setminus \{0\},
\]
where $\delta(x)\ge 0$ satisfies $\lim\limits_{x\to 0}\frac {\delta
(x)}{|x|}=0$. Then
$$
\lim\limits_{r\to 0}\inf\limits_{B_r}u>a.
$$
\label{lemma0}
\end{lem}
A slightly weaker version of Lemma \ref{lemma0} is
\begin{lem}(\cite{LL2})\ For $n\ge 2$, $R>0$, let $u\in C^2(B_R\setminus\{0\})$
satisfying $\Delta u\le 0$ in $B_R\setminus\{0\}$. Assume that there
exist $w,~v\in C^1(B_R)$ satisfying
\[
w(0)=v(0),\quad\nabla w(0)\neq\nabla v(0),
\]
and
\[
u\ge w,\quad u\ge v,\quad\mbox{in}~B_R\setminus\{0\}.
\]
Then
\[
\liminf\limits_{x\to 0}u(x)>w(0).
\]
\label{lemma00}
\end{lem}

The way we use Lemma \ref{lemma0} is as follows.  For 
some function $u$ as in the lemma, we construct
a family of $C^1$ functions $\{w^{(x)}\}$ satisfying
$$
u\ge w^{(x)}\ \ \mbox{in}\ B_1\setminus\{0\}
$$
and
$$
w^{(x)}(0)=\liminf_{|y|\to 0} u(y).
$$
An application of the lemma yields, for some $V\in \Bbb R^n$,
$$
\nabla w^{(x)}(0)=V\ \mbox{for all } x.
$$
The above could contain much information.  

To better illustrate the idea, we give
a

\noindent{\bf Proof of
 Corollary \ref{cor9}
in the case (i).}\  For $|x|$ small, shift $v$ by $x$ to obtain $v(x+\cdot)$,
which may not be $\le u$. Lower the graph of $v(x+\cdot)$  and then move it
up until one can not move further without 
cutting through the graph of $u$.  We have obtained 
$$
w^{(x)}:=v(x+\cdot)+\bar\lambda(x),
$$
which satisfies, for small $x$, 
$$
u\ge w^{(x)}\ \ \mbox{in}\ B_1\setminus\{0\}
$$
and
$$
\inf_{ B_1\setminus\{0\} }\left[ u-w^{(x)}\right]=0.
$$
By the smallness of $x$, the touching of the graphs of $u$ and
$w^{(x)}$ can not occur on $\partial B_1$.
The touching can not occur in $B_1\setminus\{0\}$ either, 
in view of the strong maximum principle.  Thus we have
$$
w^{(x)}(0)=\liminf_{ |y|\to 0} u(y).
$$
According to  Lemma \ref{lemma0}, $\nabla w^{(x)}(0)=\nabla v(x)$ is independent of $x$ and,
consequently,
 $v\equiv v(0)+\nabla v(0)\cdot x$ in $B_\epsilon$ for some
$\epsilon>0$.  Now we have $\Delta (u-v)\le 0$ and $u-v>0$
in $B_\epsilon\setminus\{0\}$, and 
(\ref{3-2y}) follows.

\vskip 5pt
\hfill $\Box$
\vskip 5pt

The paper is organized as follows.  In Section 2, we prove 
 Theorem \ref{thm1}.  In Section 3, we prove  Theorem \ref{thm7},
 Theorem \ref{thm8}, Theorem \ref{thm9}
 and Theorem \ref{thm10}.
 In Section 4, we prove Theorem \ref{thm2}.
 In Section 5, we prove  Theorem \ref{thm3} and Theorem \ref{thm11}.
 In Section 6, we prove Theorem \ref{thm4}.
 In Section 7, we prove  Theorem \ref{thmM1}, Corollary \ref{lemA-1},
 Corollary \ref{corM11}
 and  Corollary \ref{thm5}.
 In Section 8, we comment on the sharpness of  Theorem \ref{thm2}.
 In Section 9, we prove Corollary \ref{corA1}, Corollary \ref{cor9}
and Corollary \ref{cor1.6}.

Theorem \ref{thm2}, Theorem \ref{thm3}, 
Theorem \ref{thm4} and Corollary \ref{thm5}
were announced  at the international
conference in honor of Haim Brezis$'$s sixtyth birthday
 in Paris,  June 9-13, 2004.

\section{ Proof of Theorem \ref{thm1}}

\noindent{\bf Proof of Theorem \ref{thm1}.}\
 We prove
 it by contradiction.  Suppose the contrary of (\ref{3-2y}),
 then (\ref{3-2}) holds.

We first give three lemmas.
For $\epsilon>0$, let
$\lambda_\epsilon:= 1-\sqrt \epsilon.
$

\begin{lem} 
There exists some $\bar \epsilon\in (0,1)$ such that
$$
v^{x, \lambda_\epsilon}(y)<u(y),\qquad
\forall\ |x|<\epsilon\le \bar \epsilon, \
0<|y|\le 1.
$$
\label{lem5-1}
\end{lem}

\noindent{\bf Proof.}\ 
Let $\delta, \epsilon_0>0$ be some small constants chosen later,
we have,
for 
$|x|<\epsilon<\epsilon_0$ and
$0<|y|<\delta$,
\begin{eqnarray*}
v^{x, \lambda_\epsilon}(y)-u(y)
&\le &  \lambda_\epsilon^{  \frac{n-2}2 }v(x+\lambda_\epsilon y)-v(y)
\nonumber\\
&=& [1-\frac {n-2}2\sqrt\epsilon
+O(\epsilon)]
[v(y)+O(|x-\sqrt \epsilon y|]
-v(y)\nonumber\\
&=&  -\frac {n-2}2\sqrt\epsilon
v(y)+ \sqrt \epsilon  O(\sqrt \epsilon+\delta).
\label{6-1}
\end{eqnarray*}
Thus, for some small enough $\epsilon_0, \delta>0$,
\begin{equation}
v^{x, \lambda_\epsilon}(y)<u(y),\qquad
\forall\ 0<|x|<\epsilon<\epsilon_0,\
\forall\ 0<|y|<\delta.
\label{7-1}
\end{equation}

For the above $\epsilon_0$ and $\delta$,
$$
v^{x, \lambda_\epsilon}(y)=  
 v(y) +O(\sqrt \epsilon), \qquad \quad
\forall\ |x|<\epsilon<\epsilon_0, \delta\le |y|\le 1.
$$
Fix some small  $\bar \epsilon\in (0, \epsilon_0)$ so that
$$
O(\sqrt {\bar \epsilon })<
\min_{\delta \le |z|\le 1}[u(z)-v(z)].
$$
Then,
 for
$|x|<\epsilon<\bar \epsilon$ and $\delta\le |y|\le 1$,
\begin{equation}
v^{x, \lambda_\epsilon}(y)=v(y)+O(\sqrt \epsilon)<
v(y)+[u(y)-v(y)]=u(y).
\label{aaa}
\end{equation} 
Lemma \ref{lem5-1} follows from (\ref{7-1}) and (\ref{aaa}).

\vskip 5pt
\hfill $\Box$
\vskip 5pt

\begin{lem} 
There exists $\epsilon_1\in (0,1)$ such that
\begin{equation}
v^{x,\lambda}(y)<u(y),\
\
\forall\ 0<\epsilon<\epsilon_1, 1-\sqrt \epsilon\le \lambda \le 1+\sqrt \epsilon,
\ |x|<\epsilon, |y|=1.
\label{9-1}
\end{equation}
\label{lem9-2}
\end{lem}

\noindent{\bf Proof.}\
Since $v^{0,1}=v$ and
$\displaystyle{
\min_{  |y|=1 }[u(y)-v(y)]
>0 }$,
(\ref{9-1}) follows from the  continuity of $v$.

 \vskip 5pt
 \hfill $\Box$
 \vskip 5pt

\begin{lem}Under the contradiction hypothesis (\ref{3-2}),
there exists $\epsilon_2\in (0,1)$ such that
$$
\sup_{0<|y|\le 1} \left\{ v^{x, 1+\frac{\sqrt \epsilon}2}(y)
-u(y)\right\}>0,
\qquad \forall\ |x|<\epsilon<\epsilon_2.
$$
\label{lem10-3}
\end{lem}

\noindent{\bf Proof.}\
For $|x|<\epsilon<\epsilon_2$, we have, using 
(\ref{3-2}),
\begin{eqnarray*}
&&\limsup_{ |y|\to 0}
 \left\{ v^{x, 1+\frac{\sqrt \epsilon}2}(y)
-u(y)\right\}\\
&=& v^{ x, 1+\frac {\sqrt \epsilon}2 }(0)-v(0)=
\left[ 1+ \frac {n-2}2 \frac {\sqrt\epsilon} 2+O(\epsilon)\right]
[v(0)+O(\epsilon)]-v(0)\\
&=& \frac{  (n-2) \sqrt\epsilon }4 v(0)
+O(\epsilon)>0,
\end{eqnarray*}
provided that $\epsilon_2$ is small.  Lemma \ref{lem10-3}
 is established.

\vskip 5pt
\hfill $\Box$
\vskip 5pt

Now we complete the proof of Theorem \ref{thm1}.  
  Let $\bar \epsilon, \epsilon_1$ and $\epsilon_2$ be the 
  constants in Lemma \ref{lem5-1}, Lemma \ref{lem9-2}
  and Lemma \ref{lem10-3}, and let
$$
\epsilon:=\frac 18\min\{\bar \epsilon, \epsilon_1, \epsilon_2\}.
$$
For  
$|x|<\epsilon$, we know from
Lemma \ref{lem5-1} that
$$
v^{ x, 1-\sqrt\epsilon}(y)<u(y),
\qquad \forall\ 0<|y|\le 1.
$$
Thus we can define,
 for $|x|<\epsilon$,
$$
\bar \lambda(x):=\sup\{
\mu\ge 1-\sqrt \epsilon\ |\
v^{x,\lambda}(y)<u(y),
\ \forall \ 0<|y|\le 1,
\ \forall\ 1-\sqrt\epsilon \le \lambda\le \mu\}.
$$
Clearly,
\begin{equation}
\bar \lambda(x)\ge 1-\sqrt \epsilon\qquad\forall\ |x|<\epsilon.
\label{bbb}
\end{equation}
By Lemma \ref{lem10-3},
\begin{equation}
\bar \lambda (x)\le 1+\frac {\sqrt\epsilon}2,\qquad
\forall\ |x|<\epsilon.
\label{13-2}
\end{equation}
By the definition of $\bar\lambda(x)$,
\begin{equation}
v^{ x, \bar\lambda(x)}(y)\le u(y),
\qquad\forall \ |x|<\epsilon, \ \forall\
0<|y|\le 1.
\label{14-1}
\end{equation}
By Lemma \ref{lem9-2}, in view of (\ref{bbb})
 and (\ref{13-2}),
\begin{equation}
v^{x, \bar\lambda(x)}(y)<u(y),
\qquad\forall\ |x|<\epsilon, \ \forall\ 
|y|=1.
\label{14-2}
\end{equation}
By the invariance property
of $T$ and by (\ref{2-5}), 
\begin{equation}
T\left( 
v^{x, \bar\lambda(x)},
\nabla v^{x, \bar\lambda(x)},
\nabla ^2 v^{x, \bar\lambda(x)}\right)\le 0,
\qquad
\mbox{in}\ B_{\frac 32}, \ \forall\ |x|<\epsilon.
\label{14-3}
\end{equation}
In view of (\ref{14-3}), (\ref{2-5}), (\ref{14-1}) and
(\ref{14-2}), we apply the strong maximum principle
to obtain
\begin{equation}
v^{x, \bar\lambda(x)}(y)<u(y),\qquad
\forall\ |x|<\epsilon,\ \forall\
0<|y|\le 1.
\label{15-1}
\end{equation}
By (\ref{15-1}) and the definition  of $\bar \lambda(x)$, 
\begin{equation}
\liminf_{ y\to 0}
\left[ u(y)- v^{x, \bar\lambda(x)}(y)\right]
=0,\qquad\forall\ |x|<\epsilon.
\label{15-2}
\end{equation}

In view of (\ref{14-1}), (\ref{15-2}) and (\ref{2-3}),
we apply Lemma \ref{lemma00} as in \cite{LL3} 
to obtain, for some  constant vector  $V\in \Bbb R^n$,
\begin{equation}
\nabla v^{ x, \bar \lambda (x) }(0)=V,\qquad
\forall \ |x|<\epsilon.
\label{16-1}
\end{equation}

Recall 
$$
v^{ x, \bar\lambda(x)}(y)=\bar\lambda(x)^{ \frac {n-2}2 }
v(x+\bar\lambda(x)y).
$$
By (\ref{15-2}),
\begin{equation}
\alpha:=\liminf_{ y\to 0}u(y)
=v^{ x, \bar \lambda(x) }(0)=\bar \lambda(x)^{ \frac{n-2}2}
v(x),
\quad \forall\ |x|<\epsilon.
\label{16-2}
\end{equation}

So, using (\ref{16-1}) and (\ref{16-2}),
$$
V= \nabla v^{x, \bar\lambda(x)}(0)=
\bar\lambda(x)^{ \frac n2} \nabla v(x)
= \alpha^{ \frac n{n-2} }
v(x)^{ -\frac n{n-2} }\nabla v(x),
$$
i.e.
$$
\nabla\left\{\frac {n-2}2 \alpha^{ \frac n{n-2} }
v(x)^{ -\frac 2{n-2} }+V\cdot x\right\}=0,
\qquad\forall\ |x|<\epsilon.
$$
This implies, for some constant vector
$\widetilde V\in \Bbb R^n$,
$$
v(x)\equiv v(0)[1-\widetilde V\cdot x]^{ -\frac {n-2}2 },
\qquad\forall\ |x|<\epsilon.
$$
It follows that 
\begin{equation}
\Delta v(x)\ge 0, \qquad\forall\ |x|<\epsilon.
\label{hhh1}
\end{equation}
It is well known that (\ref{2-2}), (\ref{2-3}) and (\ref{hhh1})
imply
(\ref{3-2y}), contradicting to (\ref{3-2}).
Theorem \ref{thm1} is established.

\vskip 5pt
\hfill $\Box$
\vskip 5pt

\section{ Proof of Theorem \ref{thm7},
Theorem \ref{thm8}, Theorem \ref{thm9}
and Theorem \ref{thm10}}

First we give
the

\noindent{\bf Proof of Theorem \ref{thm7}.}\
The proof is similar to that  of Theorem \ref{thm1}.
Suppose the contrary of (\ref{3-2y}), then (\ref{3-2}) holds.
We still use the notation $\lambda_\epsilon:=1-\sqrt{\epsilon}$.

Instead of Lemma \ref{lem5-1} we have
\begin{lem}
There exists some small $\bar \epsilon>0$ 
 such that
 $$
 v^{x, \lambda_\epsilon}_\varphi(y)<u(y),\qquad
 \forall\ |x|<\epsilon\le \bar \epsilon, \
y\in \overline \Omega\setminus\{0\}.
 $$
 \label{lem5-1z}
 \end{lem}

 \noindent{\bf Proof.}\
Let $\delta, \epsilon_0>0$ be some small constants chosen later,
we have,
for
$|x|<\epsilon<\epsilon_0$ and
$0<|y|<\delta$,
\begin{eqnarray*}
v^{x, \lambda_\epsilon}_\varphi(y)-u(y)
&\le&   \varphi(\lambda_\epsilon)v(x+\lambda_\epsilon y)-v(y)
\nonumber\\
&=& [\varphi(1)-\varphi'(1)\sqrt{\epsilon}+\circ(1)\sqrt{\epsilon}]
[v(y)+O(|x-\sqrt \epsilon y|]
-v(y)\nonumber\\
&=& [-\varphi'(1)v(y)+\circ(1)]\sqrt{\epsilon} 
+O(\delta \sqrt{\epsilon}),
\end{eqnarray*}
where $\circ(1)\to 0$ as $\epsilon\to 0$.
Thus, for some small enough $\epsilon_0, \delta>0$,
\begin{equation}
v^{x, \lambda_\epsilon}_\varphi(y)<u(y),\qquad
\forall\ 0<|x|<\epsilon<\epsilon_0,\
\forall\ 0<|y|<\delta.
\label{7-1z}
\end{equation}

For the above $\epsilon_0$ and $\delta$,
$$
v^{x, \lambda_\epsilon}_\varphi(y)=  
\varphi(\lambda_\epsilon) v(x+\lambda_\epsilon y)
=
 v(y) +O(\sqrt \epsilon), \quad
 \forall\ |x|<\epsilon<\epsilon_0, y\in \overline \Omega
 \setminus B_\delta.
$$
Fix some small  $\bar \epsilon\in (0, \epsilon_0)$ so that
$$
O(\sqrt {\bar \epsilon })<
\inf_{\Omega\setminus B_\delta}[u(z)-v(z)].
$$
Then,
 for
 $|x|<\epsilon<\bar \epsilon$ and $y
 \in \overline \Omega\setminus B_\delta$,
 \begin{equation}
 v^{x, \lambda_\epsilon}_\varphi(y)=v(y)+O(\sqrt \epsilon)<
 v(y)+[u(y)-v(y)]=u(y).
 \label{aaaz}
 \end{equation}
 Lemma \ref{lem5-1z}  follows from (\ref{7-1z}) and (\ref{aaaz}).

\vskip 5pt
\hfill $\Box$
\vskip 5pt

\begin{lem}  Under the contradiction hypothesis
(\ref{3-2}), 
there exists $\epsilon_2\in (0,1)$ such that
$$
\sup_{y\in\Omega\setminus\{0\}}
\left\{ v^{x, 1+\frac{\sqrt \epsilon}2}_\varphi(y)
-u(y)\right\}>0,
\qquad \forall\ |x|<\epsilon<\epsilon_2.
$$
\label{lem10-3z}
\end{lem}

\noindent{\bf Proof.}\
For $|x|<\epsilon<\epsilon_2$, we have, using
(\ref{3-2}), 
\begin{eqnarray*}
&&\limsup_{ |y|\to 0}
\left\{ v^{x, 1+\frac{\sqrt \epsilon}2}_\varphi(y)
-u(y)\right\}\\
&=& v^{x, 1+\frac{\sqrt \epsilon}2}_\varphi(0)
-v(0)
= \left[
\varphi(1)+ \frac{\sqrt \epsilon}2 \varphi'(1)+\circ(\sqrt\epsilon)\right]
[v(0)+O(\epsilon)]-v(0)\\
&=& \frac{\sqrt \epsilon}2 \varphi'(1)v(0)+ \circ(\sqrt\epsilon)>0,
\end{eqnarray*}
where we have used $\varphi'(1)>0$ and
$\epsilon_2$ small.
Lemma \ref{lem10-3z} is established.

\vskip 5pt
\hfill $\Box$
\vskip 5pt

Now we complete the proof of Theorem \ref{thm7}.  
  Let
  \begin{equation}
 0< \epsilon\le \frac 18\{\bar \epsilon, \epsilon_1, \epsilon_2,
  (\epsilon_3)^2\}
  \label{71new}
  \end{equation}
  such that
  \begin{equation}
  \frac 12\le \varphi(\lambda)\le 2,\ 
  \varphi'(\lambda)>\frac 12\varphi'(1)>0,\qquad
  \forall\ |\lambda-1|\le \sqrt\epsilon.
  \label{kkk9z}
  \end{equation}
For  
$|x|<\epsilon$, we know from
Lemma \ref{lem5-1z} that
$$
v^{ x, 1-\sqrt\epsilon}_\varphi(y)<u(y),
\qquad \forall\ y\in \overline \Omega\setminus\{0\}.
$$
Thus we can define,
 for $|x|<\epsilon$,
$$
\bar \lambda(x):=\sup\{
\mu\ge 1-\sqrt \epsilon\ |\
v^{x,\lambda}_\varphi(y)<u(y),
\ \forall \ y\in \overline \Omega\setminus\{0\},
\ \forall\ 1-\sqrt\epsilon \le \lambda\le \mu\}.
$$
Clearly,
\begin{equation}
\bar \lambda(x)\ge 1-\sqrt \epsilon\qquad\forall\ |x|<\epsilon.
\label{bbbz}
\end{equation}
By Lemma \ref{lem10-3z},
\begin{equation}
\bar \lambda (x)\le 1+\frac {\sqrt\epsilon}2,\qquad
\forall\ |x|<\epsilon.
\label{13-2z}
\end{equation}
By the definition of $\bar\lambda(x)$,
\begin{equation}
\inf_{\Omega\setminus\{0\}}\left[u-
v^{ x, \bar\lambda(x)}_\varphi\right]
=0, 
\qquad\forall \ |x|<\epsilon.
\label{14-1z}
\end{equation}
By  (\ref{gg5}), in view of (\ref{14-1z}),
\begin{equation}
\liminf_{ y\to 0}
\left[ u(y)- v^{x, \bar\lambda(x)}_\varphi(y)\right]
=0,\qquad\forall\ |x|<\epsilon.
\label{15-2z}
\end{equation}

In view of (\ref{14-1z}), (\ref{15-2z}) and (\ref{2-3}),
we apply Lemma \ref{lemma0}
to obtain, for some  constant vector  $V\in \Bbb R^n$,
\begin{equation}
\nabla v^{ x, \bar \lambda (x) }_\varphi(0)=V,\qquad
\forall \ |x|<\epsilon.
\label{16-1z}
\end{equation}

Recall
$$
v^{ x, \bar\lambda(x)}_\varphi(y)=
\varphi\left(\bar\lambda(x)\right)
v(x+\bar\lambda(x)y).
$$
By (\ref{15-2z}) and (\ref{3-2}),
\begin{equation}
\alpha:=v(0)=\liminf_{ y\to 0}u(y)
=v^{ x, \bar \lambda(x) }_\varphi(0)=\varphi(\bar \lambda(x))
v(x),
\quad \forall\ |x|<\epsilon.
\label{16-2z}
\end{equation}
So, using (\ref{16-1z}) and (\ref{16-2z}),
\begin{eqnarray*}
V&=& \nabla v^{x, \bar\lambda(x)}_\varphi(0)=
\bar\lambda(x)
\varphi(\bar\lambda(x))\nabla v(x)=
\frac \alpha{v(x)}\cdot \varphi^{-1}(\frac \alpha{v(x)})
\nabla v(x).
\end{eqnarray*}

Let
$$
\psi(s):=\int_{v(0)}^s \frac \alpha t
\varphi^{-1}(\frac \alpha t)dt,
$$
we have
$$
V=\nabla_x \psi(v(x))\quad \forall\ |x|<\epsilon,
$$
i.e.
\begin{equation}
\psi(v(x))+V\cdot x=0\quad\forall\
|x|<\epsilon.
\label{hhh7}
\end{equation}

Since $\varphi$ is  $C^1$ and
$\varphi'>0$, we know that
$\psi$ is  $C^2$,
\begin{equation}
\psi'(s)= \frac \alpha s\varphi^{-1}(\frac \alpha s)>0,
\ \ \mbox{and}\ \ 
\psi''(s)=
-\frac \alpha{s^2}
\left\{
\varphi^{-1}(\frac \alpha s)+
\frac \alpha s(\varphi^{-1})'(\frac \alpha s)\right\}<0.
\label{74new}
\end{equation}
Note that we have used
(\ref{kkk9z}) in deriving the second inequality above.

Since $\psi\in C^2$ and  $\psi'>0$,
we see from (\ref{hhh7}) that $v$ is $ C^2$ near the
$0$. 
Applying $\Delta$ to (\ref{hhh7})
leads to
\begin{equation}
\psi'(v(x))\Delta v(x)+\psi''(v(x))|\nabla v(x)|^2
=0.
\label{74newnew}
\end{equation}
This implies that $\Delta v(x)\ge 0$ for  $ x$ close to
$0$.
This, together with (\ref{2-2newnew}) and (\ref{2-3newnew}), yields
(\ref{3-2y}) which contradicts to the contradiction hypothesis  
(\ref{3-2}).  Impossible.
Theorem \ref{thm7} is established.

\vskip 5pt
\hfill $\Box$
\vskip 5pt

Now we give the

\noindent{\bf Proof of Theorem \ref{thm8}.}\
The proof is similar to that  of Theorem \ref{thm7}.
We suppose that  (\ref{3-2}) holds, and we will derive
a contradiction.
We first give two lemmas whose proofs are almost identical 
to the proofs of Lemma \ref{lem5-1z} and Lemma \ref{lem10-3z}.

\begin{lem}
There exists some $\bar \epsilon>0$
 such that
$$
   v^{x, 1+\sqrt{\epsilon}}_\varphi(y)<u(y),\qquad
    \forall\ |x|<\epsilon\le \bar \epsilon, \
     y\in \overline \Omega\setminus\{0\}.
$$
        \label{lem5-1zz}
	 \end{lem}

\begin{lem}
There exists $\epsilon_2>0$ such that
$$
\sup_{y\in \Omega\setminus\{0\}}
\left\{ v^{x, 1-\frac{\sqrt \epsilon}2}_\varphi(y)
-u(y)\right\}>0,
\qquad \forall\ |x|<\epsilon<\epsilon_2.
$$
\label{lem10-3zz}
\end{lem}

Let $\epsilon$ be defined by (\ref{71new}).
For
$|x|<\epsilon$, we know from
Lemma \ref{lem5-1zz} that
$$
v^{ x, 1+\sqrt\epsilon}_\varphi(y)<u(y),
\qquad \forall\ y\in \overline \Omega\setminus\{0\}.
$$
Thus we can define,
 for $|x|<\epsilon$,
 $$
 \bar \lambda(x):=\inf\{
 \mu\ge 1+\sqrt \epsilon\ |\
 v^{x,\lambda}_\varphi(y)<u(y),
 \ \forall \ y\in \Omega\setminus\{0\},
 \ \forall\ \mu\le \lambda
 \le 1+\sqrt\epsilon\}.
 $$
 Clearly,
$$
 \bar \lambda(x)\le 1+\sqrt \epsilon\qquad\forall\ |x|<\epsilon.
$$
 By Lemma \ref{lem10-3zz},
$$
 \bar \lambda (x)\ge 1-\frac {\sqrt\epsilon}2,\qquad
 \forall\ |x|<\epsilon.
$$
 By the definition of $\bar\lambda(x)$,
$$
 v^{ x, \bar\lambda(x)}_\varphi(y)\le u(y),
 \qquad\forall \ |x|<\epsilon, \ \forall\
 y\in \overline\Omega\setminus\{0\}.
$$
The arguments between (\ref{14-1z}) and (\ref{hhh7}) 
yield (\ref{gg7}).  If $\varphi'(1)<-1$, then
$\varphi^{-1}(1)+(\varphi^{-1})'(1)=
1+\varphi'(1)^{-1}>0$,
and therefore, by (\ref{74new}), $\psi''(s)<0$ for $s$ close to
$v(0)$.  By (\ref{74newnew}), we still have
$\Delta v\ge 0$ near the origin, and we obtain (\ref{3-2y})
as usual.
Theorem \ref{thm8} is established.

\vskip 5pt
\hfill $\Box$
\vskip 5pt

\noindent{\bf Proof of Theorem \ref{thm9}.}\
Suppose the contrary of (\ref{3-2y}),   then (\ref{3-2}) holds.

\begin{lem} There exists some $\bar\epsilon>0$ such that
$$
\Phi(v,x,1-\sqrt\epsilon; y)<u(y),
\qquad\forall\ |x|<\epsilon\le \bar\epsilon, y\in \overline \Omega
\setminus\{0\}.
$$
\label{lemL4}
\end{lem}

\noindent{\bf Proof.}\ Use notation $\lambda_\epsilon=1-\sqrt\epsilon$. 
Let $\delta, \epsilon_0>0$ be some small constants chosen later,
we have,
for
$|x|<\epsilon<\epsilon_0$ and
$0<|y|<\delta$,
\begin{eqnarray*}
&&\Phi(v,x,\lambda_\epsilon; y)-u(y)\\
&\le&
 [\varphi(1)-\varphi'(1)\sqrt{\epsilon}]
[v(y)+O(\delta\sqrt\epsilon)]
-\psi'(1)\sqrt\epsilon+\circ(\sqrt\epsilon)-v(y)\\
&=& [-\varphi'(1)v(0)-\psi'(1)]\sqrt\epsilon+
\circ(\sqrt\epsilon)+O(\delta\sqrt\epsilon).
\end{eqnarray*}
Thus,  for some small enough $\epsilon_0, \delta>0$,
$$
\Phi(v,x,\lambda_\epsilon; y)<u(y),\qquad
\forall\ 0<|x|<\epsilon<\epsilon_0,\
\forall\ 0<|y|<\delta.
$$
For the above $\epsilon_0$ and $\delta$, 
$$
\Phi(v,x,\lambda_\epsilon; y)=
 v(y) +O(\sqrt \epsilon), \quad
  \forall\ |x|<\epsilon<\epsilon_0, y\in \overline \Omega
   \setminus B_\delta.
   $$
  Lemma \ref{lemL4} follows from arguments in the proof of
  Lemma \ref{lem5-1z} .

  \vskip 5pt
  \hfill $\Box$
  \vskip 5pt

\begin{lem} Under the contradiction hypothesis (\ref{3-2}), there exists
$\epsilon_2>0$ such that
$$
\sup_{  y\in \Omega\setminus\{0\}}
\bigg\{ \Phi(v,x,1+ \frac{\sqrt\epsilon}2; y)-u(y)\bigg\}>0,
\qquad \forall\ |x|<\epsilon<\epsilon_2.
$$
\label{lemL6}
\end{lem}

\noindent{\bf Proof.}\ For $|x|<\epsilon<\epsilon_2$, we have,
using (\ref{3-2}),
\begin{eqnarray*}
\limsup_{  |y|\to 0  }
\bigg\{  \Phi(v, x,1+  \frac{\sqrt\epsilon}2; y)-u(y)\bigg\}
&=&  \Phi(v, x,  1+\frac{\sqrt\epsilon}2; 0)-v(0)
\\
&=& [\varphi'(1) v(0)+
\psi'(1)]  \frac{\sqrt\epsilon}2+\circ(\sqrt\epsilon)>0,
\end{eqnarray*}
provided $\epsilon_2$ is small.

\vskip 5pt
\hfill $\Box$
\vskip 5pt

Now we complete the proof of Theorem \ref{thm9}.  Let
$$
0<\epsilon\le \frac 18\min\{\bar\epsilon, \epsilon_1, \epsilon_2,
(\epsilon_4)^2\}
$$
such that
(\ref{Q1-1}) and (\ref{Q1-2}) hold in $(1-2\epsilon, 1+2\epsilon)$.

For $|x|<\epsilon$, we know from Lemma \ref{lemL4} 
that
$$
 \Phi(v, x,1 -\sqrt \epsilon; y)<u(y),
 \qquad\forall\ y\in \Omega\setminus\{0\}.
 $$
 Thus we can define
 $$
 \bar\lambda(x):=\sup\{\mu\ge1 -\sqrt\epsilon\
 |\  \Phi(v, x, \lambda;y)<u(y),
 \forall\ y\in  \Omega\setminus\{0\},
 1-\sqrt\epsilon\le \lambda\le \mu\}.
 $$
 It follows, using Lemma \ref{lemL6}, that
 \begin{equation}
 |\bar\lambda(x)-1|\le \sqrt\epsilon,\qquad
 \forall\ |x|<\epsilon.
 \label{L7-1}
 \end{equation}
 By the definition of $\bar\lambda(x)$,
 \begin{equation}
 \inf_{ \Omega\setminus\{0\}}
 \left[ u- \Phi(v, x, \bar\lambda(x); \cdot)\right]=0.
 \label{L8-1}
 \end{equation}
By (\ref{gg5z}), in view of (\ref{L8-1}),
\begin{equation}
 \liminf_{ |y|\to 0}
\left[ u(y)-  \Phi(v, x, \bar\lambda(x); y) \right]=0,\qquad
\forall\ |x|<\epsilon.
\label{L8-3}
\end{equation}
In view of (\ref{L8-1}), (\ref{L8-3}) and (\ref{2-3newnew}), 
we obtain, using
Lemma \ref{lemma0}, that for some
constant vector $V\in \Bbb R^n$,
$$
\nabla_y  \Phi(v, x, \bar\lambda(x); y)\bigg|_{y=0}=V,\qquad
\forall\ |x|<\epsilon,
$$
i.e.
\begin{equation}
V=\varphi(\bar\lambda(x))\xi(\bar\lambda(x))\nabla v(x),
\qquad\forall\ |x|<\epsilon.
\label{Q3-2}
\end{equation}
We also know from (\ref{L8-3}) that
\begin{equation}
\liminf_{ |y|\to 0}
u(y)= \varphi(\bar\lambda(x)) v(x)+ \psi(\bar\lambda(x)),                          \qquad\forall\ |x|<\epsilon.
\label{Q3-1}
\end{equation}
Note that (\ref{3-2}) implies
$$
\liminf_{ |y|\to 0}u(y)=
v(0)=\varphi(1)v(0)+\psi(1).
$$
By (\ref{Q1-1}) and (\ref{Q3-1}), using the implicit function theorem,
$\bar\lambda(x)$ depends $C^1$  on $v$, so $\bar\lambda$ is
$C^1$, and $\bar\lambda(0)=1$. 
By (\ref{Q3-2}), we know that $\nabla v$ is $C^1$, so $v$ is $C^2$.
Applying $div$ to (\ref{Q3-2}) leads to 
\begin{equation}
0=
(\varphi\xi)(\bar\lambda(x)) \Delta v(x)
+(\varphi\xi)'(\bar\lambda(x)) \nabla \bar\lambda(x)\cdot \nabla v(x).
\label{Q5-1}
\end{equation}
Applying $\nabla$ to (\ref{Q3-1}) gives
$$
0=\varphi(\bar\lambda(x)) \nabla v(x)+
\left[ \varphi'(\bar\lambda(x)) 
v(x)+\psi'(\bar\lambda(x)) \right]\nabla \bar\lambda(x).
$$
Taking inner product of the above with $\nabla \bar\lambda(x)$, we have
\begin{equation}
0=\varphi(\bar\lambda(x)) \nabla v(x)\cdot \bar\lambda(x)+
\left[ \varphi'(\bar\lambda(x)) 
v(x)+\psi'(\bar\lambda(x)) \right]
|\nabla \bar\lambda(x)|^2.
\label{83new}
\end{equation}
This implies that $\nabla v(x)\cdot \bar\lambda(x)\le 0$ and therefore,
in view of (\ref{Q5-1}), $\Delta v(x)\ge 0$ near the origin.
This, together with $\Delta u(x)\le 0$ and $u-v>0$ for
$0<|x|<\epsilon$, yields 
 (\ref{3-2y}) violating the contradiction 
 hypothesis (\ref{3-2}).
 Impossible.  Theorem \ref{thm9} is established.

\vskip 5pt
\hfill $\Box$
\vskip 5pt

\noindent{\bf Proof of Theorem \ref{thm10}.}\  We assume that
(\ref{3-2}) holds, otherwise
we are done.
For $\epsilon>0$, let $\lambda_\epsilon:=1-\sqrt\epsilon$.

\begin{lem} There exists some $\bar \epsilon>0$ such that
$$
\Phi(v,x,\lambda_\epsilon;y)<u(y), 
\qquad
    \forall\ |x|<\epsilon\le \bar \epsilon, \
    y\in  \Omega\setminus\{0\}.
$$
	             \label{lem5-1zzz}
		              \end{lem}

\noindent{\bf Proof.}\
Since $v$ is $C^1$, we can find
small $\bar \epsilon>0$ such that
\begin{eqnarray*}
\Phi(v,x,\lambda_\epsilon;y)&=&
\left[ 1-\varphi'(1)\sqrt\epsilon\right]
\left[ v(y)-\xi'(1)\sqrt\epsilon \nabla v(y)\cdot y\right]
-\psi'(1)\sqrt\epsilon+\circ(\sqrt\epsilon)\\
&=&v(y) -\left[\varphi'(1)v(y)+\xi'(1)\nabla v(y)\cdot y
+\psi'(1)\right]\sqrt\epsilon
+\circ(\sqrt{\epsilon})\\ 
&<&u(y),\qquad \forall\ 0<|x|<\epsilon\le \bar \epsilon, y\in 
\Omega\setminus \{0\}.
\end{eqnarray*}
Lemma \ref{lem5-1zzz}
is established.

\vskip 5pt
\hfill $\Box$
\vskip 5pt

\begin{lem} Under the contradiction hypothesis (\ref{3-2}), there exists
$\epsilon_2>0$ such that
$$
\sup_{  y\in \Omega\setminus\{0\}}
\bigg\{ \Phi(v,x,1+ \frac{\sqrt\epsilon}2; y)-u(y)\bigg\}>0,
\qquad \forall\ |x|<\epsilon<\epsilon_2.
$$
\label{lemL6new}
\end{lem}

\noindent{\bf Proof.}\ The proof of Lemma  \ref{lemL6}
works here.

\vskip 5pt
\hfill $\Box$
\vskip 5pt

Follow, with obvious modification, the proof 
 of Theorem \ref{thm9} from the line after the
proof of Lemma \ref{lemL6}   until
 ``$\Delta v(x)\ge 0$ near the origin'' towards the end.
We know $\Delta u(x)\le 0$ and
$(u-v)(x)\ge 0$ for
$0<|x|<\epsilon$.  By
the mean value theorem, either $u-v>0$ on
$B_\epsilon\setminus\{0\}$ or $u-v \equiv 0$
on $B_\epsilon\setminus\{0\}$.
We know that  $u-v>0$ on
$B_\epsilon\setminus\{0\}$ 
would imply (\ref{3-2y}) and would violate
the hypothesis (\ref{3-2}), so we must have $u-v\equiv 0$ 
on $B_\epsilon\setminus\{0\}$.
Thus $\Delta v(x)=0$ on $B_\epsilon$.
With this, we deduce from (\ref{Q5-1}) and (\ref{83new}) that
$|\nabla \bar\lambda|=0$ 
in $B_\epsilon$, i.e., $\bar \lambda=\bar\lambda(0)=1$ in 
$B_\epsilon$.  Now we see from (\ref{Q3-1}) that 
$v=v(0)$ in $B_\epsilon$.
 Theorem \ref{thm10} is established.

\vskip 5pt
\hfill $\Box$
\vskip 5pt

\section{ Proof of Theorem \ref{thm2} }

\noindent{\bf Proof of Theorem \ref{thm2}.}\
Suppose the contrary of (\ref{22-4}), then there exist
some $\{x_j\}$ satisfying
$$
|x_j|\to \infty\qquad \mbox{as}\ j\to\infty,
$$
\begin{equation}
|x_j|^{ \frac {n-2}2 }u(x_j)\to \infty\qquad \mbox{as}\ j\to\infty.
\label{23-2}
\end{equation}
Consider
$$
v_j(x):=\left( 
\frac{ |x_j| }2 -|x-x_j|\right)^{ \frac{n-2}2 }u(x),\qquad
|x-x_j|\le \frac{ |x_j|}2.
$$
Let $|\bar x_j-x_j|< \frac{ |x_j| }2 $  satisfy
$$
v_j(\bar x_j)=\max_{ |x-x_j|\le \frac{ |x_j|}2  }v_j(x),
$$
and let
$$
2\sigma_j:= \frac{ |x_j| }2 - |\bar x_j-x_j|.
$$
Then
\begin{equation}
0<2\sigma_j\le \frac{ |x_j| }2.
\label{24-1}
\end{equation}
We know 
$$
(2\sigma_j)^{ \frac{n-2}2 }u(\bar x_j)
=v_j(\bar x_j)\ge v_j(x)\ge
(\sigma_j)^{ \frac{n-2}2 }u(x),\quad
\forall\ |x-\bar x_j|\le \sigma_j.
$$
Thus
\begin{equation}
u(\bar x_j)\ge 2^{ \frac{2-n}2 }u(x),
\qquad \forall\ |x-\bar x_j|\le \sigma_j.
\label{24-2}
\end{equation}
On the other hand,
 by (\ref{23-2}),
\begin{equation}
(2\sigma_j)^{ \frac {n-2}2 }u(\bar x_j)=
v_j(\bar x_j)\ge v_j(x_j)
= \left( \frac{ |x_j| }2\right)^{ \frac {n-2}2 }
u(x_j)\to \infty.
\label{25-1}
\end{equation}
Now, consider
$$
w_j(y):= \frac 1{  u(\bar x_j)  }
u\left(  \bar x_j+ \frac y{  
u(\bar x_j)^{ \frac 2{n-2}  }  }\right),
\qquad y\in \Omega_j,
$$
where
$$
\Omega_j:= \left\{
y\in \Bbb R^n\ |\ \bar x_j+
 \frac y{
 u(\bar x_j)^{ \frac 2{n-2}  }  }\ \in\
 \Bbb R^n\setminus \overline B_1\right\}.
 $$
By (\ref{24-2}) and (\ref{25-1}),
\begin{equation}
w_j(y)\le 2^{ \frac{n-2}2 },
\quad \forall\ |y|\le R_j:=\sigma_j u(\bar x_j)^{\frac 2{n-2}}
\to\infty.
\label{26-1}
\end{equation}
Since
$u(z)\ge \frac 1C>0$ for all
$|z|=1$,
we have
\begin{equation}
w_j(y)\ge \frac 1{   Cu(\bar x_j)  },
\qquad
\forall\ y\in \partial \Omega_j,
\label{26-2}
\end{equation}
where
$$
 \partial \Omega_j=\left\{ y\in \Bbb R^n\
 |\ \left|\bar x_j+ \frac y{  u(\bar x_j)^{ \frac 2{n-2} }  }\right|=1
 \right\}.
 $$
For any $y\in \partial \Omega_j$,
\begin{equation}
|\frac y{  u(\bar x_j)^{ \frac 2{n-2} } } |
\ge |-\bar x_j|-
|\bar x_j + \frac y{  u(\bar x_j)^{ \frac 2{n-2} } } |
=|\bar x_j|-1\ge \frac 12 |\bar x_j|.
\label{26-3}
\end{equation}
Thus, using (\ref{26-2}) and (\ref{26-3}),
\begin{eqnarray}
\min_{ y\in \partial \Omega_j} |y|^{n-2}w_j(y)&\ge &
\min_{ y\in \partial \Omega_j}  \left\{
\frac 12 |\bar x_j| u(\bar x_j)^{ \frac 2{n-2} }\right\}^{n-2}
w_j(y)\nonumber\\
&\ge&  \min_{ y\in \partial \Omega_j} 
 \left\{
 (\frac 12 |\bar x_j|)^{n-2} u(\bar x_j)^2\right\}
 \frac 1{ Cu(\bar x_j) }\nonumber\\
& =& (\frac 12)^{n-2} \frac 1C |\bar x_j|^{n-2} u(\bar x_j).
 \label{27-1}
 \end{eqnarray}
Clearly
\begin{equation}
|\bar x_j|\ge \frac 12 |x_j|.
\label{28-0}
\end{equation}
We deduce from the above and (\ref{24-1}) that
\begin{equation}
|\bar x_j|\ge 2\sigma_j.
\label{28-00}
\end{equation}
Thus
\begin{eqnarray}
|\bar x_j|^{n-2} u(\bar x_j)&\ge &
(\frac 12 |x_j|)^{ \frac{n-2}2 }
|\bar x_j|^{ \frac{n-2}2 }u(\bar x_j)
\ge (\frac 12 |x_j|)^{ \frac{n-2}2 }
(2\sigma_j)^{ \frac{n-2}2 } u(\bar x_j)
\nonumber\\
&=& |x_j|^{ \frac{n-2}2 }
(R_j)^{ \frac{n-2}2 } \to\infty.
\label{28-1}
\end{eqnarray}
We deduce from (\ref{27-1}) and (\ref{28-1}) that
\begin{equation}
\lim_{j\to\infty} \min_{ y\in \partial \Omega_j}
|y|^{n-2} w_j(y)=\infty.
\label{28-2}
\end{equation}
By (\ref{26-3}) and (\ref{28-00}),
\begin{equation}
|y|\ge  \frac 12 |\bar x_j| u(\bar x_j)^{ \frac 2{n-2} }
\ge \sigma_j u(\bar x_j)^{ \frac 2{n-2} }
= R_j, \qquad \forall\ y\in \partial \Omega_j.
\label{29-1}
\end{equation}
By (\ref{22-1}),
 (\ref{22-3}) and
the invariance of the equation, and by (\ref{26-1})
and (\ref{28-2}),
\begin{equation}
\left\{
\begin{array}{l}
F(A^{w_j})=1, \ \  \ A^{w_j}\in U, w_j>0, \Delta w_j\le 0\ \mbox{in}\ \Omega_j,\\
w_j(0)=1,\\
w_j(y)\le 2^{  \frac {n-2}2 }, \ \ \ 
\forall\ |y|\le R_j,\\
\min_{ y\in \partial \Omega_j} \left\{
|y|^{n-2} w_j(y) \right\}\to\infty.
\end{array}
\right.
\label{30-1}
\end{equation}

For all $|x|< \frac{R_j}{ 10}$,
let
$$
(w_j)_{x,\lambda}(y):=
\left(  \frac \lambda { |y-x| }\right)^{n-2}
w_j\left(  x+  \frac{\lambda^2(y-x) }{  |y-x|^2 }\right)
$$
and 
$$
\bar \lambda_j(x)=\sup\{
\mu>0\ |\ 
(w_j)_{x,\lambda}(y)\le 
w_j(y),\ \forall\ |y-x|\ge \lambda,
y\in \overline \Omega_j, \forall\ 0<\lambda<\mu\}>0,
$$
is well defined, see proof of lemma 2.1 in \cite{LZ}.

For  $|x|< \frac{R_j}{ 10}$, $0<\lambda\le \frac{R_j}4$
and $y\in \partial \Omega_j$, we know,
from  (\ref{29-1}) that 
$$|y-x|\ge |y|-|x|\ge \frac {9}{10}|y|\ge \frac 9{10}R_j,
$$
and
$$
\bigg| x+\frac { \lambda^2(y-x) }{  |y-x|^2 }\bigg|\le
|x|+\frac { \lambda^2}  {|y-x|}
\le \frac{R_j}{10}+
(\frac {R_j}4)^2 (\frac{10}{9 R_j})\le
\frac {R_j}2.
$$

So
$$
|y|^{n-2} (w_j)_{x,\lambda}(y)
\le  2^{ \frac{n-2}2}  (\frac{10}{9})^{n-2} \lambda^{n-2},
\qquad \mbox{for}\  |x|<\frac {R_j}{10},
0<\lambda\le \frac{R_j}4, y\in \partial\Omega_j.
$$
Because of the last line in  (\ref{30-1}), there exist $r_j\to \infty$,
$r_j\le \frac{ R_j }4$, such that
\begin{equation}
 (w_j)_{x,\lambda}<w_j\qquad \mbox{on}\ 
 \partial \Omega_j\ \mbox{for all}\ |\lambda|\le r_j.
\label{104new}
\end{equation} 
Namely, for all $ |\lambda|\le r_j$, no touching
of $(w_j)_{x,\lambda}$ and $w_j$ can occur on $\partial \Omega_j$.

Now we prove
\begin{equation}
\bar\lambda_j(x)\ge r_j.
\label{L10-1}
\end{equation}
Suppose the contrary, $\bar\lambda_j<r_j$,
then, in view of (\ref{104new}), we
can use the strong maximum principle and
the Hopf Lemma as in the proof of lemma 2.1 in \cite{LL1} to
show
\begin{equation}
(w_j)_{  x, \bar\lambda(x)}<w_j\qquad
\mbox{in}\ \Omega_j\setminus
\overline {  B_{\bar\lambda(x)}(x)  },
\label{L10-2}
\end{equation}
\begin{equation}
\frac {\partial }{
\partial \nu}
\left[  w_j- 
(w_j)_{  x, \bar\lambda(x)}
\right]
\bigg|_{  \partial B_{ \bar\lambda(x) }(x)  }
>0,
\label{L10-3}
\end{equation}
where $\frac \partial {\partial \nu}$ denotes differentiation in outer 
normal direction of 
$B_{ \bar\lambda(x) }(x)$. 
Applying Theorem \ref{thm0} to the Kelvin transformation of $w_j$ and
$(w_j)_{ x, \bar(x) }$ which turn
the singularity of $w_j$ from $\infty$ to
$0$, we have
\begin{equation}
\liminf_{  |y|\to \infty }
|y|^{n-2} \left[
w_j(y)- (w_j)_{ x, \bar\lambda(x)}(y)\right]>0.
\label{L11-1}
\end{equation}
As usual, (\ref{104new}), (\ref{L10-2}),
(\ref{L10-3}) and (\ref{L11-1})
allow the moving sphere procedure
to go beyond $\bar\lambda(x)$, contradicting to
the definition of $\bar\lambda(x)$.
We have established (\ref{L10-1}).  Once we have (\ref{L10-1}),
the
  argument in the proof of theorem 1.2 in
\cite{LL3} then leads to contradiction.
Theorem \ref{thm2} is established.

\vskip 5pt
\hfill $\Box$
\vskip 5pt

\section{ Proof of Theorem \ref{thm3} and Theorem \ref{thm11}}

\noindent{\bf Proof of Theorem \ref{thm3}.}\
By the positivity and the superharmonicity of
$u$ in $\Bbb R^n\setminus\{0\}$,
$$
\liminf_{|y|\to 0}u(y)>0,
\qquad \liminf_{ |y|\to \infty}
|y|^{n-2} u(y)>0.
$$
For all $|x|>0$, 
 we can
prove as usual, see e.g. \cite{LZ} or \cite{LL1},
that there exists
$\lambda_0(x)\in (0, |x|)$ such that
for all $0<\lambda <\lambda_0(x)$,
$$
u_{ x, \lambda}(y)
:=\left(  \frac\lambda{|y-x|}\right)^{n-2}
u\left(x+\frac{\lambda^2(y-x)}{  |y-x|^2 }\right)\le u(y),\ \
\forall\ |y-x|\ge \lambda, |y|>0.
$$
Define
$$
\bar \lambda(x)=\sup
\{ 0<\mu<|x|\ |
\ u_{x,\lambda}(y)\le u(y),
\forall\ |y-x|\ge \lambda, |y|\ne 0,
0<\lambda<\mu\}.
$$
We will prove 
\begin{equation}
\bar\lambda(x)=|x|, \qquad \forall\ |x|>0.
\label{105new}
\end{equation}
Suppose for some $|x|>0$,
$\bar \lambda(x)<|x|$, then we obtain,  using
  the strong maximum principle and the Hopf Lemma as
in section 2 of \cite{LL1} 
 and in view of (\ref{34-1}), 
\begin{equation}
u(y)-u_{x, \bar\lambda(x)}(y)>0, \quad
\forall\ |y-x|>\bar \lambda(x), |y|\ne 0,
\label{mm1}
\end{equation}
and 
\begin{equation}
\partial_\nu [u-u_{x, \bar\lambda(x)}]
\bigg|_{ \partial B_{\bar \lambda(x)} }
>0,
\label{mm2}
\end{equation}
where $\partial_\nu$ denotes the unit outer normal
derivative.

By Theorem \ref{thm0} with $v=u_{x,\lambda}$, 
\begin{equation}
\liminf_{ |y|\to 0}
[ u(y)- u_{x,\lambda}(y)]>0.
\label{mm4}
\end{equation}

Applying Theorem \ref{thm0} with $u(y)$ replaced by
$|y|^{2-n}u(\frac y{ |y|^2})$ and $v(y)$ by
$|y|^{2-n}u_{x,\lambda}(\frac y{ |y|^2})$ leads to
\begin{equation}
\liminf_{ |y|\to \infty} 
\left( |y|^{n-2} [ u(y)- u_{x,\lambda}(y)] \right)>0.
\label{mm3}
\end{equation}

But this would violate the definition of $\bar\lambda(x)$, since 
(\ref{mm1}), (\ref{mm2}), (\ref{mm3}) and (\ref{mm4})
would allow the moving sphere procedure to continue beyond
$\bar \lambda(x)$.
Thus we have proved (\ref{105new}).

It follows that
\begin{equation}
u_{x, |x|}(y)\le u(y),\qquad \forall\ |y-x|\ge |x|>0.
\label{mm5}
\end{equation}
For any unit vector $e\in \Bbb R^n$,  for any $y\in \Bbb R^n$ satisfying
$y\cdot e<0$, and for any $R>0$, we have, by
(\ref{mm5}) with $x=Re$,
$$
u(y)\ge u_{x, |x|}(y)
= \left( \frac {|x|} { |y-x| }\right)^{n-2}
u\left( x+  \frac{  |x|^2 (y-x) }
{|y-x|^2 } \right).
$$
Sending $R$ to infinity in the above leads to
$$
u(y)\ge u(y-2 (y\cdot e)e).
$$
This gives the radial symmetry of the function $u$.
Theorem \ref{thm3} is established.

\vskip 5pt
\hfill $\Box$
\vskip 5pt

\noindent{\bf Proof of Theorem \ref{thm11}.}\
As usual, 
$$
\liminf_{|y|\to 0}u(y)>0
$$
and, for all $0<|x|<\frac 12$,
$$
\bar \lambda(x)=\sup
\{ 0<\mu<|x|\ |
\ u_{x,\lambda}(y)\le u(y),
\forall\ |y-x|\ge \lambda, 0<|y|\le 1,
0<\lambda<\mu\}>0
$$
is well defined.

For $|y|=1$, and $0<\lambda<|x|<\frac 12$,
$$
\bigg| \bigg\{
x+\frac{ \lambda^2(y-x) }{ |y-x|^2 } \bigg\}-x\bigg|\le 4\lambda^2\le 4|x|^2.
$$
So
$$
\bigg| \bigg\{
x+\frac{ \lambda^2(y-x) }{ |y-x|^2 } \bigg\}-x\bigg|\le 
\frac {|x|}4,\quad 
\forall\ 0< \lambda<|x|<\frac 14.
$$
Thus, by Theorem \ref{thm2}$'$, 
$$
u\left(x+\frac{ \lambda^2(y-x) }{ |y-x|^2 }\right)
\le C|x|^{\frac{2-n}2},
$$
and, for some $\epsilon>0$,  
$$
u_{x,\lambda}(y)\le
C\lambda^{n-2}|x|^{\frac{2-n}2}
\le C|x|^{ \frac{n-2}2 }<u(y),
\qquad\forall \ 0<\lambda<|x|\le \epsilon,
|y|=1.
$$
This means that no touching of $u_{x,\lambda}$ and $u$ may
occur on $\partial B_1$ in the moving sphere procedure. By the strong
maximum principle as usual,  the moving sphere procedure
can not stop due to touching of  $u_{x,\lambda}$ and $u$ 
in $B_1\setminus\{0\}$.
On the other hand, by 
Theorem \ref{thm0},  
no touching of  $u_{x,\lambda}$ and $u$  at the origin may occur.
Therefore, $\bar\lambda(x)=|x|$ for all $|x|\le \epsilon$.
We have proved (\ref{radial}). 
Let
$\displaystyle{
v(y):=|y|^{2-n}u(\frac y{|y|^2})}$,  (\ref{radial})
amounts to the following:
$$
v(y)\le v(y_\lambda)\quad \forall\ y\cdot e\ge \frac 1\epsilon,
 e\in \Bbb R^n, |e|=1,
$$
where $y_\lambda=y+2(\lambda-x\cdot e)e$ is the 
reflection of $y$ in the plane $x\cdot \tau =\lambda$.
Now we can follow the proof of corollary 6.2 in \cite{CGS} to
obtain
 (\ref{vv9}).
Theorem \ref{thm11} is established.

\vskip 5pt
\hfill $\Box$
\vskip 5pt

\section{  Proof of Theorem \ref{thm4}  }

\noindent{\bf Proof of Theorem \ref{thm4}.}\
By (\ref{41-1}) and the fact that
$A^u\in U$, $\Delta u\le 0$ in
$B_2\setminus\{0\}$. By Theorem \ref{thm2}$'$,
\begin{equation}
\sup_{ 0< |x|\le 1} |x|^{ \frac {n-2}2 }u(x)<\infty.
\label{42-1}
\end{equation}
Since
 $\Delta u\le 0$ in
 $B_2\setminus\{0\}$, we have
 \begin{equation}
 u(x)\ge \min_{\partial B_1} u>0,
 \qquad \forall\ 0<|x|\le 1.
 \label{42-2}
 \end{equation}
 Let
 $$
 \xi(x)=\frac{n-2}2 u(x)^{ -\frac 2{n-2} },
 \qquad 0<|x|<1,
 $$
 we have, as in the proof of lemma 6.5 in \cite{LL1},
 \begin{equation}
 (D^2\xi)\ge \frac 1{n-2} u^{  \frac { 2-2n}{n-2}  }
 |\nabla u|^2 I,\qquad
 B_1\setminus\{0\}.
 \label{43-1}
 \end{equation}
 We know from (\ref{42-2}) that
 $$
 \xi(x)\le C, \qquad \mbox{on}\ 0<|x|\le 1.
 $$
 Here and throughout the rest of the proof of  Theorem \ref{thm4},
 $C>1$ denotes some positive constant which
 may change its value from line to line. The constant
 $C$ is allowed  depend on $u$.

 By the convexity of $\xi$  --- see (\ref{43-1}),
 \begin{equation}
 |\nabla \xi(x)|\le C,\qquad
 \forall\ 0<|x|\le \frac 12,
 \label{43-2}
 \end{equation}
 and $\xi$ can be extended as a Lipschitz function in $B_{\frac 12}$.

 Clearly $0\le \xi\le C$ on $B_{\frac 12}$.

 We divide into two cases:

 $\underline{  \mbox{Case 1} }.$\ $\xi(0)>0$,

  $\underline{  \mbox{Case 2} }.$\ $\xi(0)=0$.

  In Case 1,
  \begin{equation}
  0<\frac 1C \le \xi
  <C<\infty\quad\mbox{on}\ B_{\frac 12}.
  \label{44-1}
  \end{equation}
  By (\ref{44-1}) and (\ref{43-2}),
  $$
  \frac 1C\le u\le C\quad
  \mbox{and}\quad |\nabla u|\le C,\qquad \mbox{on}\ B_{\frac 12},
  $$
We arrive at the conclusion of Theorem \ref{thm4}.

We need to rule out the possibility of Case 2.  In Case 2,
we have, by (\ref{43-2}),
$$
0<\xi(x)\le C|x|,\qquad \forall\ 0<|x|<\frac 12,
$$
i.e.
$$
u(x)\ge \frac 1C |x|^{ -\frac{n-2}2 },\qquad
\forall\ 0<|x|<\frac 12.
$$
This and (\ref{42-1}) give
\begin{equation}
\frac 1C |x|^{ -\frac{n-2}2 }
\le u(x)\le C|x|^{ -\frac{n-2}2 },\qquad
\forall\ 0<|x|<\frac 12.
\label{45-2}
\end{equation}
Since
$$
u=(\frac 2{n-2})^{ -\frac{n-2}2 }\xi^{ -\frac {n-2}2 },
$$
we have, for some constant $a>0$,
$$
u^{ \frac{2-2n}{n-2} }|\nabla u|^2
=(n-2)a \xi ^{n-1} (\xi^{ -\frac n2} |\nabla \xi|)^2
=(n-2)a\xi^{-1}|\nabla \xi|^2.
$$
Thus, by (\ref{43-1}),
\begin{equation}
(D^2\xi)\ge a\xi^{-1}
|\nabla \xi|^2 I,\qquad
\mbox{in}\ B_{\frac 12}\setminus\{0\}.
\label{46-1}
\end{equation}
Fixing $e=(1, 0, \cdots, 0)$, and let
$$
f(t)=\xi(te), \qquad 0<t<\frac 12.
$$
Then
$$
f'(t)=\xi_1(te), \quad f''(t)=\xi_{11}(te),
$$
and, by (\ref{46-1}),
\begin{equation}
f''(t)\ge a \xi(te)^{-1}
|\nabla \xi(te)|^2
\ge a\xi(te)^{-1} |\xi_1(te)|^2
=a f(t)^{-1} f'(t)^2,
\quad 0<t<\frac 12.
\label{46-2}
\end{equation}
$\underline{ \mbox{Claim} }.$\ $f'(t)>0$, $\forall\ 
0<t<\frac 12$.

$\underline{ \mbox{Proof of the Claim} }.$\ For all $0<t<\frac 12$, there
exists some $0<s<t$ such that
\begin{equation}
f'(s)=\frac{ f(t)-f(0) }{ t-0}
=\frac{  f(t)  }t>0.
\label{47-1}
\end{equation}
By (\ref{46-2}), $f''\ge 0$ on
$(0, \frac 12)$.  So, since $s<t$, we have
\begin{equation}
f'(s)\le f'(t).
\label{47-2}
\end{equation}
The above claim follows from (\ref{47-2}) and
(\ref{47-1}).

Because of the Claim, we 
rewrite (\ref{46-2}) as
$$
\frac{ f''}{ f'}
\ge \frac{ a f'}{ f}\qquad
\mbox{on}\
(0, \frac 12),
$$
or
$$
(\log f')'\ge (a\log f)'\qquad
\mbox{on}\ (0,\frac 12).
$$
For any $0<s<t <\frac 12$,
we deduce from the above that
\begin{equation}
\log f'(t)-\log f'(s)\ge
a[ \log f(t)-\log f(s)].
\label{48-1}
\end{equation}
By (\ref{45-2}),
$$
\frac \tau C\le f(\tau)
\le C\tau, \qquad \forall\ 0<\tau <\frac 12.
$$
For all $0<\tau<s$, there exists
some $0<\theta =\theta(\tau)<\tau$ such that
\begin{equation}
\frac 1C \le \frac{ f(\tau) }\tau
=\frac{ f(\tau) -f(0) }{ \tau -0}
=f'(\theta(\tau)).
\label{49-2}
\end{equation}
Since $\theta(\tau)<\tau<s<\frac 12$, and
since
$f''\ge 0$ on $(0, \frac 12)$, we have
\begin{equation}
f'(\theta(\tau))\le f'(s).
\label{49-3}
\end{equation}
Putting together (\ref{49-2}) and (\ref{49-3}), we have
\begin{equation}
\frac 1C \le f'(s), \qquad
\forall\ 0<s<\frac 12.
\label{50-1}
\end{equation}
By (\ref{48-1}) and (\ref{50-1}), for any
$0<s<t<\frac 12$, we have
$$
\log f(t)-\log f(s)
\le  \frac 1a \{ \log f'(t)-\log f'(s) \}
\le \frac 1a \log f'(t)+\frac 1a \log C.
$$
Fixing $t=\frac 14$ in the above, we have
$$
\log f(s)\ge -C, \qquad \forall\ 0<s<\frac 14,
$$
i.e.
$$
f(s)\ge e^{ -C}, \qquad 
\forall\ 0<s<\frac 14.
$$
Sending $s$ to $0$ leads to
$$
0=f(0)\ge e^{-C},
$$
impossible.  We have ruled out the possibility
of Case 2, and 
therefore have established Theorem \ref{thm4}.

\vskip 5pt
\hfill $\Box$
\vskip 5pt

\section{Proof of Theorem \ref{thmM1}, Corollary \ref{lemA-1},
Corollary \ref{corM11}
and  Corollary \ref{thm5}}

\noindent{\bf Proof of Theorem \ref{thmM1}.}\  The proof makes use of
arguments in the proof of theorem 2.7 in \cite{TW}.
We mainly treat the possible singularity of $\xi$ at the origin. 
We first assume in addition that $\xi\in C^2(B_2\setminus\{0\})$ and
$D^2\xi\in \overline U$ in $B_2\setminus\{0\}$.
In view of (\ref{M1-0}), we may assume that 
$D^2\xi\in U$ in $B_2\setminus\{0\}$ since
we may replace $\xi$ by 
$\xi(x)+\epsilon |x|^2$ for  $\epsilon>0$ and then send $\epsilon$ to $0$.
(\ref{M2-3}) follows from subharmonicity of $\xi$ in
$B_2\setminus\{0\}$ and the fact that $\sup_{ B_2\setminus\{0\} }\xi
<\infty$.
It is easy to see that 
we may  assume without loss of generality that 
$$
1\le \xi\le 2\quad \mbox{in}\ B_2\setminus\{0\},
$$
Fix some  $C>1$ such that for  all $0<|\bar x|<\frac 14$, 
\begin{equation}
\xi(\bar x)+C\eta(x-\bar x)>\xi(x),
\qquad \forall\ |x-\bar x|=1.
\label{A2-1}
\end{equation}
Here and throughout the proof we use $C$ to denote some constant
depending only on $\eta$ which may vary from line to line.
Consider, for $A\ge 0$,
$$
\eta_A(x):= \xi(\bar x)+C\eta(x-\bar x)+A,
\qquad |x-\bar x|\le 1.
$$
Clearly
$$
\eta_0(\bar x)=\xi(\bar x), \ \ 
\eta_2(x)\ge \xi(x),\qquad \forall\ |x-\bar x|\le 1, x\ne 0.
$$
It is easy to see that for some 
 $0\le \overline A\le 2$,
\begin{equation}
\eta_{\overline A}(x)\ge \xi(x), \qquad |x-\bar x|\le 1, x\ne 0,
\label{A2-2}
\end{equation}
and
\begin{equation}
\inf_{ |x-\bar x|\le 1, x\ne 0}
[\eta_{\overline A}(x)-\xi(x)]=0.
\label{A2-3}
\end{equation}
We must have
\begin{equation}
\eta_{\overline A}(x)> \xi(x), \qquad |x-\bar x|\le 1, x\ne 0,
x\ne \bar x.
\label{nn8}
\end{equation}
Indeed, by  (\ref{A2-1}), $\eta_{\overline A}(x)>\xi(x)$
for all $|x-\bar x|=1$.  If for some $\hat x\ne 0, \hat x\ne \bar x$ and
$|\hat x-\bar x|<1$, $\eta_{\overline A}(\hat x)=\xi(\hat x)$, then, in
view of (\ref{A2-2}), $
D^2\eta_{\overline A}(\hat x)\ge D^2 \xi(\hat x)\in U 
$ which implies, in view of (\ref{M1-0}),  $D^2\eta_{\overline A}(\hat x)\in
U$, violating (\ref{M1-2}).

We know that 
$$
\eta_{\overline A}(x)> \xi(x), \qquad 0<|x|<|\bar x|.
$$
Let $\varphi(\lambda):=\lambda$,
$u=C-\xi$, $v=C-\eta_{\bar A}$, $a=\frac 12|\bar x|$, $\Omega=B_a$, where
$C$ is some constant satisfying $C>\eta_{\bar A}$ 
in $\overline B_a$.

\noindent {\it Claim.}\ There exists $\epsilon_3>0$ such that
(\ref{gg5}) holds for the above.

\noindent{\it Proof of the Claim.}\
Suppose the contrary, then for some $\epsilon_3>0$ small and
for some $|x|<\epsilon_3$ and
 $|\lambda-1|<\epsilon_3$, 
we have
$$
\min_{\partial B_a}
[u-v^{x,\lambda}_\varphi]>
\inf_{  B_a\setminus\{0\} }
[u-v^{x,\lambda}_\varphi]
=0
$$
and
$$
\liminf_{|y|\to 0}
[u-v^{x,\lambda}_\varphi](y)>0.
$$
Then for some $0<|\bar y|<a$,
$$
[u-v^{x,\lambda}_\varphi]
(\bar y)=0.
$$
It follows that
$$
D^2u(\bar y)\ge D^2 v^{x,\lambda}_\varphi(\bar y)
$$
i.e.
$$
\lambda^2 \varphi(\lambda) D^2 \eta_{ \overline A}(x+\lambda \bar y)\ge D^2 \xi(\bar y)
\in U.
$$
It follows, using (\ref{M1-0}) and
(\ref{M1-0new}), $D^2\eta_{ \overline A}(x+\lambda \bar y)\in U$,
contradicting to (\ref{M1-2}).
The Claim has been proved.

Now we apply
 Theorem \ref{thm7}
to obtain
  $$
  \liminf_{ |x|\to 0}
    [\eta_{\overline A}(x)-\xi(x)]>0.
 $$
Thus, using also (\ref{A2-3}), (\ref{A2-2}) and (\ref{nn8}),
we have
$$
\eta_{\overline A}(\bar x)=\xi(\bar x),
$$
i.e. $\overline A=0$.  Then by (\ref{A2-2}),
\begin{equation}
\xi(x)\le \xi(\bar x)+C\eta(x-\bar x),\qquad
\forall\ |x-\bar x|\le 1, x\ne 0.
\label{131}
\end{equation}
Since (\ref{131}) holds for all 
$0<|\bar x|, |x|<\frac 14$, switching the roles of
$\bar x$ and $x$, we obtain
$$
|\xi(x)-\xi(\bar x)|
\le C[\eta(x-\bar x)+\eta(\bar x-x)],
\qquad \forall\ |x-\bar x|\le 1,
\qquad
\forall\ 0<|\bar x|, |x|<\frac 14.
$$

Now we complete the proof of Theorem \ref{thmM1}:
(\ref{M2-3}) still follows from (\ref{M1-3}) and the fact
that $\xi$ is bounded from above.  Let $\{\xi_i\}$ be in
$C^2(B_2\setminus\{0\})$ such that
(\ref{M2-0}), (\ref{M2-1}) and (\ref{M2-2}) hold.
We have proved (\ref{M2-4}) for $\{\xi_i\}$, with
constant $C(\eta)$ independent of $i$.  Sending
$i$ to $\infty$, we obtain 
(\ref{M2-4}) for $\xi$.  
Theorem \ref{thmM1} is established.

\vskip 5pt
\hfill $\Box$
\vskip 5pt

\noindent{\bf Proof of Corollary \ref{corM11}.}\
(\ref{M12-1}) follows from the
superharmonicity and the positivity of $u$ in
$B_2\setminus\{0\}$.  It is easy
to see that (\ref{M12-2}) implies either
(\ref{M13-1}) or (\ref{M13-2}).  By a limit procedure, as in the proof
of Theorem \ref{thmM1}, we 
only need to establish (\ref{M12-2}) for the $u_i$.  Now we 
drop the index $i$ in the notation. Let
$\xi=u^{ -\frac 2{n-2} }$, then $\xi\in L^\infty(B_2\setminus\{0\})$, 
$$
\Delta \xi\ge 0,\ D^2\xi \in \overline U,\qquad
\mbox{in}\ B_2\setminus\{0\}.
$$
Estimate (\ref{M12-2}) follows from Theorem \ref{thmM1}.

\vskip 5pt
\hfill $\Box$
\vskip 5pt

\noindent{\bf Proof of Corollary \ref{thm5}.}\
Let $U=U_k$ and 
 $\eta(x)=|x|^\alpha$.  Then it is known that 
 $\eta$ satisfies the properties in Corollary \ref{corM11}.
 Corollary \ref{thm5} follows from Corollary \ref{corM11}.

\vskip 5pt
\hfill $\Box$
\vskip 5pt

\noindent{\bf Proof of Corollary \ref{lemA-1}.}\
Let  $U=U_k$ and 
 $\eta(x)=|x|^\alpha$. It  follows from Theorem \ref{thmM1}.

\vskip 5pt
\hfill $\Box$
\vskip 5pt

\section{Sharpness of Theorem \ref{thm2}}

The two lemmas in this section give the sharpness of
 Theorem \ref{thm2} as stated in Remark \ref{rem1.3}.

\begin{lem}
For $n\ge 3$, let
$$
u(x)=|x|^{ \frac {2-n}2 }, \qquad
x\in \Bbb R^n\setminus\{0\}.
$$
Then
\begin{equation}
\lambda(A^u)\equiv \left\{
-\frac 12, \frac 12, \cdots, \frac 12\right\},
\qquad \mbox{on}\ \Bbb R^n\setminus\{0\}.
\label{37-1}
\end{equation}
\label{lem14}
\end{lem}

\noindent{\bf Proof.}\  We write $u(x)$ as $u(r)$ with $r=|x|$.
We only need to verify (\ref{37-1}) at
$x=(r, 0, \cdots, 0), r>0$.
At the point, we have, as
in the proof of theorem 1.6 in \cite{LL3},
\[
\nabla u(x)=(u'(r),0\cdots,0),\quad\nabla^2
u(x)=diag(u''(r),\frac{u'(r)}{r},\cdots, \frac{u'(r)}{r}),
\]
and
\[
A^u(x)=diag(\lambda_1^u(r),\lambda_2^u(r),\cdots,\lambda_n^u(r)),
\]
where
$$
\left\{
\begin{array}{lcl}
&&\lambda_1^u(r)=-\frac{2}{n-2}u^{-\frac{n+2}{n-2}}u''
+\frac{2(n-1)}{(n-2)^2}u^{-\frac{2n}{n-2}}(u')^2\\
&&\lambda_2^u(r)=\cdots=\lambda_n^u(r)=
-\frac{2}{n-2}u^{-\frac{n+2}{n-2}}\frac{u'}{r}
-\frac{2}{(n-2)^2}u^{-\frac{2n}{n-2}}(u')^2.
\end{array}
\right.
$$
With this we compute:

$$
u'=\frac{2-n}2 r^{-\frac n2}=\frac{2-n}2 u^{ \frac n{n-2} },
\quad
u''= -\frac n2 u^{ \frac 2{n-2} }u'=\frac{  n(n-2) }2 u^{ \frac{n+2}{n-2} }.
$$

$$
\lambda_1^u(r)= -\frac n2+ \frac{n-1}2=-\frac 12,
$$
$$
\lambda_2^u(r) =\cdots= \lambda_n(r)=
-\frac 2{n-2} u^{ -\frac {n+2}{n-2} }
(\frac {2-n}2) r^{ -\frac{n+2}2}
-\frac 2{  (n-2)^2 }(\frac {n-2}2)^2=\frac 12.
$$
Lemma \ref{lem14}
is established.

\vskip 5pt
\hfill $\Box$
\vskip 5pt

\begin{lem} For $\bar \lambda=(-1, 1, \cdots, 1)\in \Bbb R^n$, $n\ge 2$, 
\begin{equation}
\left\{
\begin{array}{rl}
\sigma_k(\bar\lambda)>0,& \mbox{for}\ 
1\le k<\frac n2,\\
\sigma_k(\bar \lambda)=0,&  \mbox{for}\
k=\frac n2,\\
\sigma_k(\bar \lambda)<0,&  \mbox{for}\
\frac n2<k\le n.
\end{array}
\right.
\label{39-1}
\end{equation}
It follows that
$
(-\frac 12, \frac 12, \cdots, \frac 12)$ belongs to $
\Gamma_k, \forall\ 1\le k<\frac n2,
$
and $
(-\frac 12, \frac 12, \cdots, \frac12) 
$ does not belong to $ \Gamma_k,  \forall\
k\ge \frac n2.
$
\label{lem39-4}
\end{lem}

\noindent{\bf Proof of Lemma \ref{lem39-4}.}\
For $n=2$ or for $k\in \{1, n\}$, (\ref{39-1}) is obvious.  In the rest
of the proof, we assume that $n\ge 3$ and $ 2\le k\le n-1$.
For $\lambda=(\lambda_1, \cdots, \lambda_n)\in \Bbb R^n$,
$$
\det\left(tI+diag(\lambda_1, \cdots, \lambda_n)\right)
=t^n +\sigma_1(\lambda)t^{n-1}
+\sigma_2(\lambda)t^{n-2}+\cdots+
\sigma_{n-1} t+\sigma_n(\lambda).
$$
Taking $\lambda=\bar \lambda$ and setting
$$
f(t):=(t-1)(t+1)^{n-1}
\equiv t^n+\sigma_1(\bar \lambda)t^{n-1}
+\cdots+
\sigma_{n-1}(\bar \lambda)t
+\sigma_n(\bar \lambda).
$$
Then
$$
\frac{d^k}{dt^k}f(0)= k! \sigma_{n-k}(\bar \lambda),
\qquad 1\le k\le n.
$$
Rewriting 
$$
f(t)=(t-1)(t+1)^{n-1}=(t+1-2)(t+1)^{n-1}
=(t+1)^n -2 (t+1)^{n-1}.
$$
Since
$$
\frac{d^k}{dt^k} (t+1)^n\bigg|_{t=0}=
n(n-1) \cdots (n-k+1),
$$
$$
\frac{d^k}{dt^k}  (t+1)^{ n-1} \bigg|_{t=0}=
(n-1)(n-2) \cdots (n-k+1)(n-k),
$$
we have
\begin{eqnarray*}
\frac{d^k}{dt^k}f(0)&=& n\{ (n-1)(n-2) \cdots (n-k+1)\}\nonumber\\
&&
-2 \{ (n-1)(n-2)\cdots (n-k+1)\}(n-k)\nonumber\\
&=&\{ (n-1)(n-2) \cdots (n-k+1)\}(2k-n).
\label{40-2}
\end{eqnarray*}
Since $ (n-1)(n-2) \cdots (n-k+1)>0$,  Lemma \ref{lem39-4} follows from the above.

\vskip 5pt
\hfill $\Box$
\vskip 5pt

\section{Proof of Corollary \ref{corA1}, Corollary \ref{cor9}
and Corollary \ref{cor1.6}}

We first give the 

\noindent{\bf Proof of Corollary \ref{cor1.6}.}\
Since
$$
\Phi(v,x,\lambda;y)=\Phi(v,0,\lambda; y+\xi(\lambda)^{-1}x),
$$
it is easy to see from (\ref{BB1}) and (\ref{2-2}) that
for some small $\epsilon_4>0$,
$$
T\left(\Phi(v,x,\lambda;\cdot), \nabla \Phi(v,x,\lambda;\cdot),
\nabla^2 \Phi(v,x,\lambda;\cdot)\right)\le 0,\ \ 
\mbox{in}\ B_{\epsilon/2},\
\forall\ |x|<\epsilon_4, |\lambda-1|<\epsilon_4,
$$
and
$$
u>\Phi(v,x,\lambda;\cdot),\ \ \mbox{on}\
\partial B_{\epsilon/2},\
\forall\ |x|<\epsilon_4, |\lambda-1|<\epsilon_4.
$$
Since
$$
T(u, \nabla u, \nabla^2 u)\ge 0,\ \ \mbox{in}\
B_{\epsilon/2}\setminus\{0\},
$$
and since the operator is elliptic, we can easily verify
(\ref{gg5z}), with $\Omega=B_{\epsilon/2}$,
by a contradiction argument using the maximum principle
on $B_{\epsilon/2}\setminus B_\delta$ for some small
$\delta>0$.  An application of Theorem \ref{thm9} yields
(\ref{3-2y}).

\vskip 5pt
\hfill $\Box$
\vskip 5pt

Now we give the 

\noindent{\bf Proof of Corollary \ref{cor9}.}\ We only need to verify that
operators $T$ satisfy the hypotheses of Corollary \ref{cor1.6}.

If $T$ satisfies (i), we let $\varphi(\lambda)\equiv \xi(\lambda)\equiv 1$ and
$\psi(\lambda)=\lambda-1$.  Then
$$
\Phi(v,0,\lambda;y)=v(y)+\lambda-1,
$$
and
\begin{eqnarray*}
&&T\left(\Phi(v,0,\lambda;\cdot), \nabla \Phi(v,0,\lambda;\cdot), \nabla^2 \Phi(v,0,\lambda;\cdot)\right)
\\
&=& S\left(\nabla \Phi(v,0,\lambda;\cdot), \nabla^2 \Phi(v,0,\lambda;\cdot)\right)
=S(\nabla v, \nabla^2 v)
= T(v,\nabla v, \nabla^2 v)\le 0.
\end{eqnarray*}
The hypotheses of Corollary \ref{cor1.6} are satisfied.

If $T$ satisfies (ii), we let $\varphi(\lambda)=\lambda$,
$\xi(\lambda)\equiv 1$ and $\psi(\lambda)\equiv 0$.  Then
$$
\Phi(v,0,\lambda;y)=\lambda v(y),
$$
and therefore, for $|\lambda-1|<\epsilon$, 
$$
sign\ T\left(\Phi(v,0,\lambda;\cdot), \nabla \Phi(v,0,\lambda;\cdot), \nabla^2 \Phi(v,0,\lambda;\cdot)\right)
=sign\ T(v,\nabla v, \nabla^2 v).
$$
The hypotheses of Corollary \ref{cor1.6} are satisfied.

If $T$ satisfies (iii), we let $\varphi(\lambda)=\lambda$, $\xi(\lambda)\equiv 1$ and
$\psi(\lambda)=\lambda-1$.  Then
$$
\Phi(v,0,\lambda;y)=\lambda v(y)+\lambda-1
$$
and
\begin{eqnarray*}
&& T\left(\Phi(v,0,\lambda;\cdot), \nabla \Phi(v,0,\lambda;\cdot), \nabla^2 \Phi(v,0,\lambda;\cdot)\right)
\\
&=& S\left( \frac 1 {\lambda v+\lambda}\cdot \lambda \nabla v,
\frac 1{\lambda v+\lambda}\cdot \lambda \nabla^2 v\right)
=S\left( \frac 1{v+1}\nabla v, \frac 1{v+1}\nabla^2v\right)
\\
&=&T(v,\nabla v, \nabla^2 v)\le 0.
\end{eqnarray*}
The hypotheses of Corollary \ref{cor1.6} are satisfied.

\vskip 5pt
\hfill $\Box$
\vskip 5pt

Before proving Corollary \ref{corA1}, we give a lemma.
For $\beta\in \Bbb R$, let
$$
\varphi_\beta(\lambda):=\lambda^\beta,\quad
v^\lambda_{ \varphi_\beta }(y)
:=\varphi_\beta (\lambda) v(\lambda y)=\lambda^\beta
 v(\lambda y).
 $$
\begin{lem}
For $n\ge 1$ and $\beta\in \Bbb R\setminus\{0\}$, 
let $T\in C^0( \Bbb R_+\times \Bbb R^n\times {\cal S}^{ n\times n})$.
Then
\begin{equation}
T(v^{\lambda}_{\varphi_\beta}, \nabla v^{\lambda}_{\varphi_\beta},
\nabla^2 v^{\lambda}_{\varphi_\beta})(\cdot)
\equiv T(v, \nabla v, \nabla^2 v)(\lambda \cdot)
\qquad \mbox{in}\ \Bbb R^n
\label{1-1w}
\end{equation}
holds for  any positive function $v\in C^2(\Bbb R^n)$ and for
any  $\lambda>0$ if and only if
\begin{equation}
T(t,p,M)\equiv S(t^{ -\frac  {1+\beta}\beta }p,
t^{ -\frac{2+\beta}\beta }M), \qquad\forall\
 (t,p,M)\in
  \Bbb R_+\times \Bbb R^n\times {\cal S}^{ n\times n}
  \label{kkk6}
  \end{equation}
  for some $S\in C^0(\Bbb R^n\times {\cal S}^{ n\times n})$.
\label{thmT}
\end{lem}

\noindent{\bf Proof.}\ Assuming 
(\ref{1-1w}), then
for any positive $C^2$ function $v$ and for all
$\lambda>0$, we know from (\ref{1-1}) that
$$
T\left( \lambda^{\beta}v(\lambda y), 
\lambda^{ 1+\beta }\nabla v(\lambda y), \lambda^{ 2+\beta}\nabla^2v
(\lambda y)\right)\equiv T\left(v(\lambda y), \nabla v(\lambda y),
\nabla^2 v(\lambda y)\right),
$$
i.e.
\begin{equation}
T(ts, t^{\frac {1+\beta}\beta}p, t^{\frac{2+\beta}\beta}M)
=T(s,p,M)\qquad
\forall\ (t,s,p,M)\in \Bbb R_+\times
\Bbb R_+\times \Bbb R^n\times {\cal S}^{ n\times n}.
\label{kkk1}
\end{equation}
Taking $t=\frac 1s$ in the above leads to
(\ref{kkk6}), with $S(p,M):=T(1,p,M)$.

On the other hand, if (\ref{kkk6}) holds for some $S$, then
\begin{eqnarray*}
T(ts, t^{\frac {1+\beta}\beta }p,
t^{  \frac {2+\beta}\beta }M)
&=&S\left( (ts)^{  -\frac  {1+\beta}\beta  }(t^{\frac 
 {1+\beta}\beta }p),
 (ts)^{  -\frac  {2+\beta}\beta  }  (t^{\frac  {2+\beta}\beta  }M\right)\\
&=& S\left( s^{  -\frac {1+\beta}\beta 
}p, s^{  -\frac  {2+\beta}\beta   }M\right)
=T(s,p,M).
\end{eqnarray*}
This implies (\ref{1-1}).
Lemma \ref{thmT} is established.

\vskip 5pt
\hfill $\Box$
\vskip 5pt

Now the

\noindent{\bf Proof of Corollary \ref{corA1}.}\
Let
$$
\varphi(\lambda):=\lambda^\beta,\quad
v_{\varphi}^{x,\lambda}(y):=
\varphi(\lambda)v(x+\lambda y),\quad
\Omega:=B_1.
$$
For $\epsilon_3>0$ small, we have, for
any $|x|<\epsilon_3$ and $|\lambda-1|<\epsilon_3$,
$$
v_\varphi^{x,\lambda}<u,\qquad\mbox{on}\ \partial B_1
$$
and, by Lemma \ref{thmT}, 
$$
T\left(v_\varphi^{x,\lambda}, \nabla v_\varphi^{x,\lambda}, \nabla^2 v_\varphi^{x,\lambda}\right)
(\cdot)\equiv T\left(v, \nabla v, \nabla^2 v\right)(x+\lambda\cdot).
$$Thus (\ref{gg5}) can be proved by a contradiction 
argument using the maximum principle since
$$
T(u, \nabla u, \nabla^2 u)
\ge 0\ge 
T\left(v_\varphi^{x,\lambda}, \nabla v_\varphi^{x,\lambda}, \nabla^2 v_\varphi^{x,\lambda}\right),
\quad \mbox{in}\ \Omega\setminus\{0\}.
$$

If $0<\beta<\infty$, then $\varphi'(1)>0$, and
we can apply Theorem \ref{thm7} to obtain (\ref{3-2y}).
If $-\infty<\beta <-1$, then $\varphi'(1)<-1$, and an application of
 Theorem \ref{thm8}
yields (\ref{3-2y}).  If $\beta=-1$, then, by 
Theorem \ref{thm8}, either (\ref{3-2y}) holds or,
for some $V\in \Bbb R^n$ and $\epsilon>0$,
$$
v(x)-v(0)+V\cdot x\equiv 0, \qquad |x|<\epsilon.
$$
The latter implies that  $\Delta(u-v)\le 0$ in $B_\epsilon\setminus \{0\}$, and
(\ref{3-2y}) follows as usual since $u-v>0$ in $B_\epsilon\setminus\{0\}$.

If $-1<\beta<0$, then by 
Theorem \ref{thm8}, either (\ref{3-2y}) holds or,
for some $V\in \Bbb R^n$ and $\epsilon>0$,
\begin{equation}
\psi(v(x))+V\cdot x=
-\beta v(0)\left[\left(\frac {v(x)}{v(0)}\right)^{ -\frac 1\beta} -1\right]+V\cdot x
\equiv 0,\qquad |x|<\epsilon.
\label{CCC1}
\end{equation}
We deduce from (\ref{CCC1}) that, in $B_\epsilon$,
$$
v(x)\equiv v(0)\left[1+
\frac{V\cdot x}{ \beta v(0)}\right]^{-\beta},
\qquad v^{ -\frac{1+\beta}\beta } \nabla v\equiv - v(0)^{ -\frac{1+\beta}\beta }V,
$$
and
$$
v^{ -\frac{2+\beta}\beta }\nabla^2 v\equiv (1+\frac 1\beta)
v(0)^{ -\frac {2+2\beta}\beta } V\otimes V.
$$
It follow, using also (\ref{2-5}) and (\ref{CCC2}), that
\begin{eqnarray*}
0&\ge & T(v, \nabla v, \nabla^2 v)= S\left(v^{ -\frac{1+\beta}\beta } \nabla v,
v^{ -\frac{2+\beta}\beta }\nabla^2 v\right)
\\
&=& S\left( v(0)^{ -\frac{1+\beta}\beta } V,
(1+\frac 1\beta) v(0)^{ -\frac {2+2\beta}\beta } V\otimes V\right)\\
&=& S\left( v(0)^{ -\frac{1+\beta}\beta } V, 0\right)+
\int_0^1 \left[\frac d{dt}
S\left( v(0)^{ -\frac{1+\beta}\beta } V,
t(1+\frac 1\beta) v(0)^{ -\frac {2+2\beta}\beta } V\otimes V\right)\right]
dt\\
&\ge &  (1+\frac 1\beta) v(0)^{ -\frac {2+2\beta}\beta }
\int_0^1 \left[\frac 
{\partial S}{\partial M_{ij}}\left( v(0)^{ -\frac{1+\beta}\beta } V,
t(1+\frac 1\beta) v(0)^{ -\frac {2+2\beta}\beta } V\otimes V\right)V_iV_j\right]
dt.
\end{eqnarray*}
Since $1+\frac 1\beta<0$ and
$
\displaystyle{
\left( -\frac{\partial S}{\partial M_{ij}}\right)>0
},$ we see from above that $V=0$, i.e. $v\equiv v(0)$, in
$B_\epsilon$.  We obtain (\ref{3-2y}) as usual.

\vskip 5pt
\hfill $\Box$


\vskip 5pt



\begin{thebibliography}{99}
\bibitem{CGS} L. Caffarelli, B. Gidas and J. Spruck,
Asymptotic symmetry
and local behavior of semilinear elliptic equations with
critical Sobolev growth, Comm. Pure Appl. Math.
 42 (1989),  271-297.
 \bibitem{CNS}  L. Caffarelli,  L. Nirenberg
 and  J. Spruck,  The Dirichlet problem for nonlinear second-order
 elliptic equations, III:  Functions of the eigenvalues of the Hessian.
 Acta Math. 155 (1985),
 261-301.
 \bibitem{C} S.Y.A. Chang, Non-linear
 elliptic equations in conformal Geometry,
 Nachdiplom Lectures  Course Notes, ETH, Zurich,
  Springer-Birkh\"auser, 2004.
   \bibitem{CGY}  S.Y.A. Chang, M. Gursky and P. Yang,
    An equation of Monge-Ampere type in conformal geometry,
     and four-manifolds of positive Ricci curvature,
      Ann. of Math. 155 (2002),  709-787.
        \bibitem{CGY2}  S.Y.A. Chang, M. Gursky and P. Yang,
	  A prior estimate for a class of nonlinear equations
	    on 4-manifolds, Journal D'Analyse Journal Mathematique
	     87 (2002),  151-186. 
	     \bibitem{CHY}  S.Y.A. Chang, Z.C. Han and P. Yang,
Classification of singular
radial solutions to the $\sigma_k-$Yamabe equation on
annular domains,  preprint.
\bibitem{E} G.C. Evans,
Potentials and positively infinite singularities of
harmonic functions, Monatsh. Math. 43 (1936),
419-424.
	     \bibitem{d1} Mar\'ia del Mar Gonz\'alez,
	     Classification of singularities for a subcritical
	     fully non-linear problem, preprint.
	      \bibitem{d2} Mar\'ia del Mar Gonz\'alez,
	     Removability of singularies
	     for a class of fully   non-linear elliptic
	     equations, preprint.
	      \bibitem{GNN}  B. Gidas, W.M. Ni and L. Nirenberg,
	       Symmetry and related properties via
	        the maximum principle, Comm. Math. Phys. 68 (1979), 209-243.
		 \bibitem{GV}  M. Gursky and J. Viaclovsky,
		  Convexity and
		   singularities of curvature equations in conformal geometry,
		   arXiv:math.DG/0504066 v1 4 Apr 2005.
	\bibitem{H} Z.C. Han, Local pointwise estimates
	for solutions of the $\sigma_2$ curvature equation
	on $4-$manifolds,  International Mathematics Research Notices
	79 (2004), 4269-4292.
	\bibitem{LL1}
		   A. Li and Y.Y. Li,
		   On some conformally invariant fully
		   nonlinear equations, Comm. Pure Appl. Math. 56 (2003), 1416-1464.
		   \bibitem{LL2}  A. Li and Y.Y. Li,
		   A general Liouville type theorem  for some conformally invariant fully
		   nonlinear equations,
		   arXiv:math.AP/0301239 v1 21 Jan 2003.
		   \bibitem{LL4} A. Li and Y.Y. Li,
		   Further results on Liouville type theorems
		    for some conformally invariant fully
		     nonlinear equations,
		      arXiv:math.AP/0301254 v1 22 Jan 2003.
		      \bibitem{LL3}  A. Li and Y.Y. Li,
		       On some conformally invariant fully nonlinear equations, Part II:
		        Liouville, Harnack and Yamabe,
			 arXiv:math.AP/0403442 v1 25 Mar 2004.
			  \bibitem{Li1} Y.Y. Li,
			  Degenerate
			  conformally invariant fully nonlinear equations,
			  preprint.
			     \bibitem{LZ} Y.Y. Li and L. Zhang,
			      Liouville type theorems and Harnack type inequalities for semilinear
			       elliptic equations, Journal d'Analyse Mathematique 90 (2003), 27-87.

 \bibitem{TW} N.S. Trudinger and X.-J. Wang,
   Hessian measures. II. Ann. of Math.
      150 (1999), 579-164.
      \bibitem{V1} J. Viaclovsky,  Estimates and existence
			       results for some fully nonlinear elliptic
			       equations on Riemannian manifolds, Comm. Anal. Geom. 10 (2002), 815-846.
\end{thebibliography}
\end{document}